"\pdfoutput=1"
\documentclass[11pt]{amsart}

\usepackage{fullpage}
\usepackage{amssymb}
\usepackage{amsmath}
\usepackage{amsxtra}
\usepackage{amscd}
\usepackage{graphicx}
\usepackage{amsfonts}
\usepackage{pb-diagram}
\usepackage{float}
\usepackage{enumerate}

\numberwithin{equation}{section}

\newtheorem{defn}{Definition}

\newtheorem{rem}{Remark}
\newtheorem{note}{Note}

\begin{document}

\title[Topological Surgery and its Dynamics]{Topological Surgery and its Dynamics}

\author{Sofia Lambropoulou}
\address{Department of Mathematics,
National Technical University of Athens,
Zografou campus, GR--157 80 Athens, Greece.}
\email{sofia@math.ntua.gr}
\urladdr{http://www.math.ntua.gr/$\sim$sofia}

\author{Stathis Antoniou}
\address{Department of Mathematics,
National Technical University of Athens,
Zografou campus, GR--157 80 Athens, Greece.}
\email{stathis.antoniou@gmail.com}

\author{Nikola Samardzija}
\address{Emerson Electric Co.,
11533 Park Ridge Dr. W Minnetonka, MN 55305, USA.}
\email{Nick.Samardzija@Emerson.net}

\thanks{This research has been co--financed by the European Union (European
Social Fund -- ESF) and Greek national funds through the Operational Program
``Education and Lifelong Learning" of the National Strategic Reference
Framework (NSRF) -- Research Funding Program: THALIS}



\date{}

\begin{abstract}
Topological surgery occurs in natural phenomena where two points are selected and attracting or repelling forces are applied.  The two points are connected via an invisible `thread'. In order to model topologically such phenomena  we introduce dynamics in 1-, 2- and 3-dimensional topological surgery, by means of attracting or repelling forces between two selected points in the manifold, and we address examples. We also introduce the notions of solid 1- and 2-dimensional topological surgery, and of truncated 1-, 2- and 3-dimensional  topological surgery, which are more appropriate for modelling natural processes. On the theoretical level, these new notions allow to visualize 3-dimensional surgery and to connect surgeries in different dimensions. 
 We hope that through this study, topology and dynamics of many natural phenomena as well as  topological surgery may now be better understood.
\end{abstract}

\maketitle

\section*{Introduction}

The aim of this study is to draw a connection between topological surgery in dimensions 1, 2 and 3 and many natural phenomena.  For this we introduce new theoretical concepts which allow to explain the topology of such phenomena via surgery and also to connect topological surgeries in different dimensions. The new concepts are the introduction of forces, attracting or repelling, in the process of  surgery, the notion of solid 1- and 2-dimensional surgery and the notion of truncated 1-, 2- and 3-dimensional surgery.

\smallbreak
Topological surgery is a technique used for changing the homeomorphism type of a topological manifold, thus for creating new manifolds out of known ones. A homeomorphism between two $n$-manifolds is a continuous bijection, such that the inverse map is also continuous. 
 Further, manifolds with homeomorphic boundary may be attached together and a homeomorphism between their boundaries can be used as `glue'.  An {\it $n$-dimensional topological surgery} on an $n$-manifold $M$ is, roughly, the topological procedure whereby an appropriate $n$-manifold with boundary is removed from  $M$  and  is replaced  by another $n$-manifold with the same boundary, using a `gluing'  homeomorphism, thus creating a new  $n$-manifold $\chi(M)$ (not necessarily different from the starting one). For details see, for example, \cite{PS, Ro}. 

\smallbreak

Apart from just being a formal topological procedure, topological surgery appears in nature in numerous, diverse processes of various scales for ensuring new results. Such processes are initiated by attracting or repelling forces between two  points, or {\it `poles'},  which seem to be joined by some {\it invisible `thread'}.  To list some examples, 1-dimensional surgery happens in DNA recombination and in the reconnection of cosmic magnetic lines. 2-dimensional surgery is exhibited in the formation of whirls, in blowing bubbles, in the Falaco solitons and in the cell mitosis. 3-dimensional surgery can be observed, for example, in the formation of tornadoes, or the magnetic field excited by a current loop. 

Surgery  in nature is usually performed on basic manifolds with boundary.  
In each dimension the basic closed (compact without boundary), connected, oriented (c.c.o.) $n$-manifold, on which surgery is usually performed, is the $n$-sphere, $S^n$, which may be viewed as ${\mathbb R}^n$ with all points at infinity compactified to one single point. We also need to recall that the basic connected, oriented $n$-manifold with boundary is the solid $n$-ball, $D^n$.  In particular for $n=3$, other 3-manifolds with boundary that we will be using are: the solid torus, which can be described as the product set $S^1\times D^2$, and the handlebodies, which generalize the solid torus, having higher genus.

We are particularly interested in situations related to  2- and 3-dimensional topological surgery of the attracting type. Here, a {\it `drilling' process} along the invisible thread seems to be initiated, resulting in passage from spherical to toroidal shape.  `Drilling' with coiling seems to be a natural choice in various physical processes, probably for being the most effective way for opening up a hole. 

\smallbreak

 From the above, topological surgery is not just a mathematical technique used for changing the homeomorphism type of a manifold. It can also serve as a mathematical tool for explaining the change of topology in many natural phenomena.  For this reason we introduce dynamics in the process of surgery.

\smallbreak

 In Sections~1, 2 and 3 we recall first the mathematical definitions of topological surgery in dimensions 1, 2 and 3 respectively.  Then, we introduce  dynamics in topological surgery distinguishing two different types: via {\it attracting forces} and via {\it repelling forces} between two selected points, the `poles'.  Each one of these two types of dynamics can be eventually viewed as the reverse of the other.  We also introduce the notions of {\it solid 1- and 2-dimensional surgery}, whereby the interior space is now filled in. Also, the notions of {\it truncated 1-, 2- and 3-dimensional surgery}, whereby surgery is being localized. All these notions are better adapted to natural or physical processes exhibiting topological surgery and we address briefly in each dimension such examples. 

\smallbreak

Moreover, truncated 3-dimensional surgery allows to visualize 3-dimensional surgery, for which the fourth dimension is otherwise needed. This is presented in Section~3.

\smallbreak

Finally, in Section~4  we use all the new notions introduced in the previous sections in order to pinn down the relation among topological surgeries in dimensions 1, 2 and 3. 

\smallbreak

 The first author was always fascinated by 3-dimensional surgery and was trying to find ways to visualize it. So Figure~\ref{3Dtruncated} dates back several years ago. Further, our work  is inspired by our connection of 3-dimensional topological surgery with a dynamical system \cite{SaGr1,SaGr2,Sa,La,SS,SSN}. Then, on one hand we will have a mathematical model for 3-dimensional surgery. On the other hand, through our connection many natural phenomena can be modelled through our dynamical system. We hope that our observations, new definitions and ideas will serve in turn as inspiration for many more interesting connections.

\section{1-dimensional topological surgery}\label{1D}

\subsection{}
\begin{figure}[ht!]
\begin{center}
\includegraphics[width=12cm]{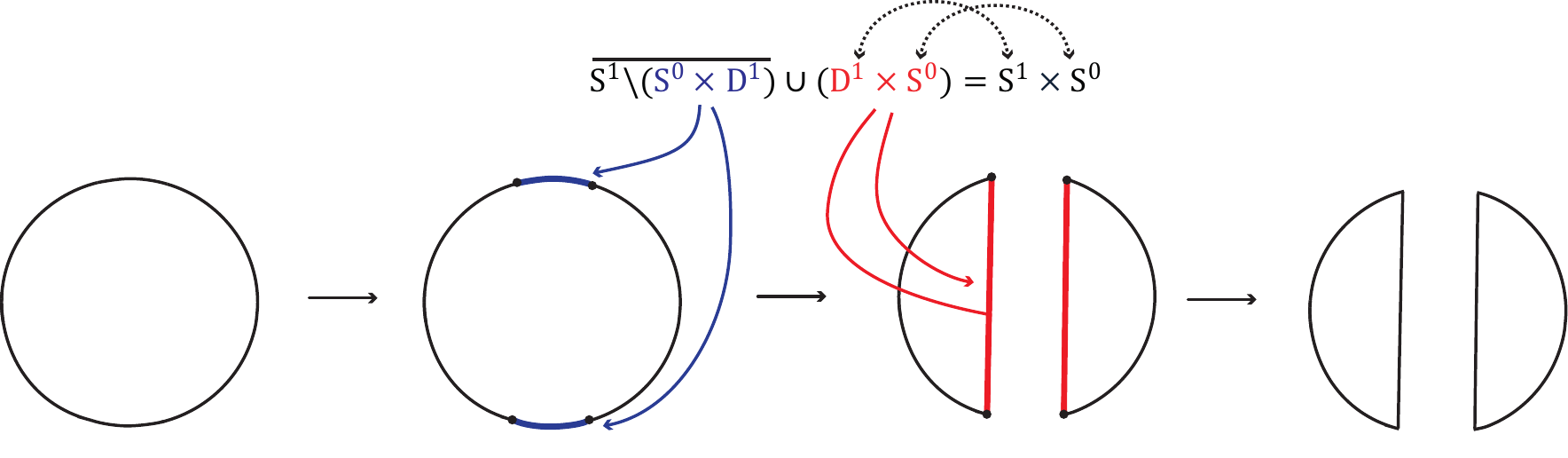}
\caption{1-dimensional surgery.}
\label{Formal1D}
\end{center}
\end{figure}

Starting with $S^1$, {\it 1-dimensional surgery} means that: two segments $S^0\times D^1$  are removed from $S^1$ and they are replaced (in the closure of the remaining manifold) by two different segments $D^1 \times S^0$ by reconnecting the four boundary points $S^0\times S^0$ in a different way. In the end we obtain two circles $S^1 \times S^0$ or one, depending on the type of reconnection, see Figure~\ref{Formal1D}. Recall that $S^0$ consists in two points.

\smallbreak

\begin{figure}[ht!]
\begin{center}
\includegraphics[width=12cm]{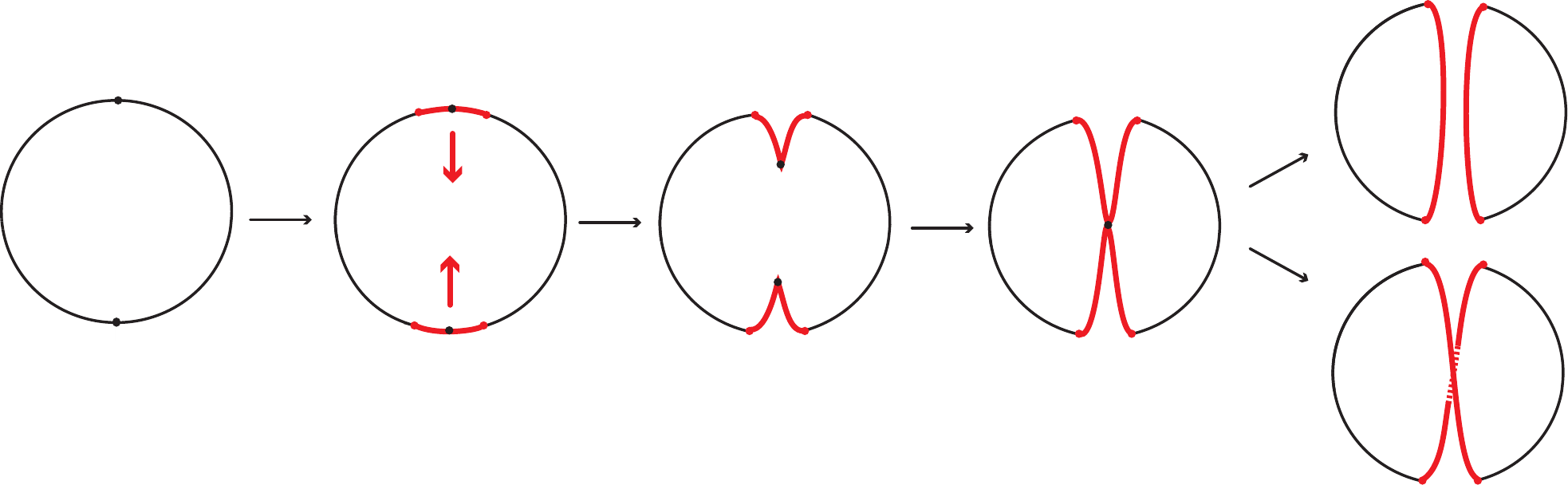}
\caption{Attracting 1-dimensional surgery.}
\label{1Dattract}
\end{center}
\end{figure}

\begin{figure}[ht!]
\begin{center}
\includegraphics[width=12cm]{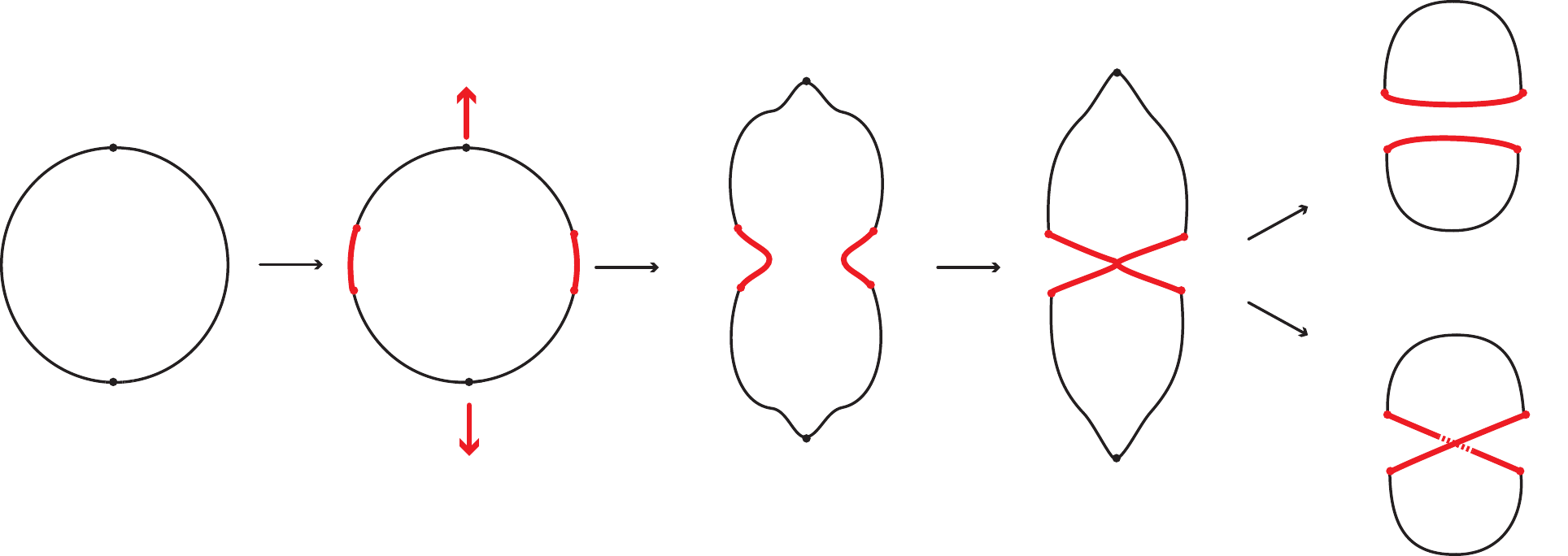}
\caption{Repelling 1-dimensional surgery.}
\label{1DRepel}
\end{center}
\end{figure}

The above definition of 1-dimensional surgery gives only the initial and the final stage. 
 In order to adress natural phenomena exhibiting 1-dimensional surgery or to understand how  1-dimensional surgery happens we need a non-static description. We will describe the process by  introducing dynamics. The process starts with two points specified on the circle, on which attracting forces are applied. Then, the two segments $S^0\times D^1$, which are neighbourhoods of the two points, get close to one another. When the two segments touch, recoupling takes place giving rise to the two final segments $D^1 \times S^0$, which split apart. See Figure~\ref{1Dattract}. This type of 1-dimensional surgery shall be called {\it attracting 1-dimensional surgery}. We also have the {\it repelling 1-dimensional surgery}, whereby  repelling forces are applied on the two points, as illustrated in Figure~\ref{1DRepel}. Note here that the recoupling does not take place between the neighbourhoods of the two repelling points but between `complementary' segments, which get closer by passive reaction.

\smallbreak

\smallbreak

\subsection{}
\begin{figure}[ht!]
\begin{center}
\includegraphics[width=11cm]{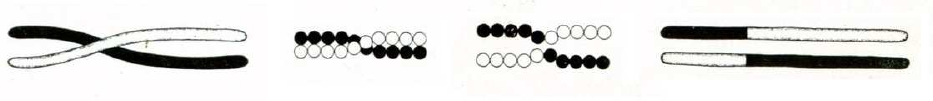}
\caption{Crossing over of chromosomes in DNA recombination.}
{\footnotesize \underline{{\it Source}}: Wikipedia}
\label{DNArecomb}
\end{center}
\end{figure}

\begin{figure}[ht!]
\begin{center}
\includegraphics[width=10cm]{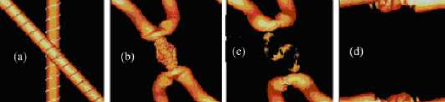}
\caption{The reconnection of cosmic magnetic lines.}
{\footnotesize \underline{{\it Source}}: R.B. Dahlburg, S.K. Antiochos, {\em Reconnection of Antiparallel Magnetic Flux Tubes}, J.  Geophysical Research {\bf 100}, No. A9 (1995) 16991--16998.}
\label{magneticrecon}
\end{center}
\end{figure}

\begin{figure}[ht!]
\begin{center}
\includegraphics[width=12cm]{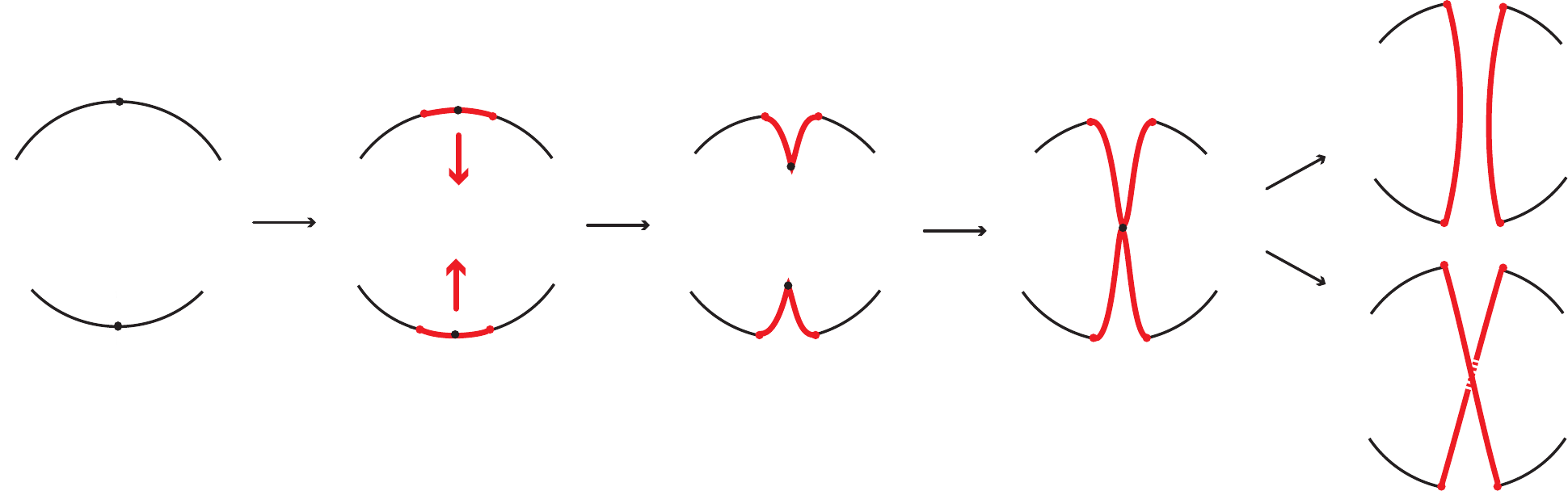}
\caption{Truncated 1-dimensional surgery by attraction.}
\label{1Dattrtrunc}
\end{center}
\end{figure}

In practice, 1-dimensional surgery happens on arcs or lines. That is, the initial space is the closure of $S^1 \setminus  (D^1 \times S^0) = S^0 \times D^1$, and we remove from it a smaller $S^0 \times D^1$. We shall call this type of surgery  {\it truncated 1-dimensional surgery}. See Figure~\ref{1Dattrtrunc} for the case of attracting forces. Truncated 1-dimensional surgery happens, for example, on the double helix and recombines DNA, thus changing the genetic sequence. See Figure~\ref{DNArecomb}.  Also, in magnetic reconnection --the phenomenon whereby cosmic magnetic field lines from different magnetic domains are spliced to one another-- changing the patterns of connectivity with respect to the sources. See Figure~\ref{magneticrecon} (cf. \cite{DaAn}).    
 
\smallbreak

\smallbreak

\smallbreak

\subsection{}
\begin{figure}[ht!]
\begin{center}
\includegraphics[width=11cm]{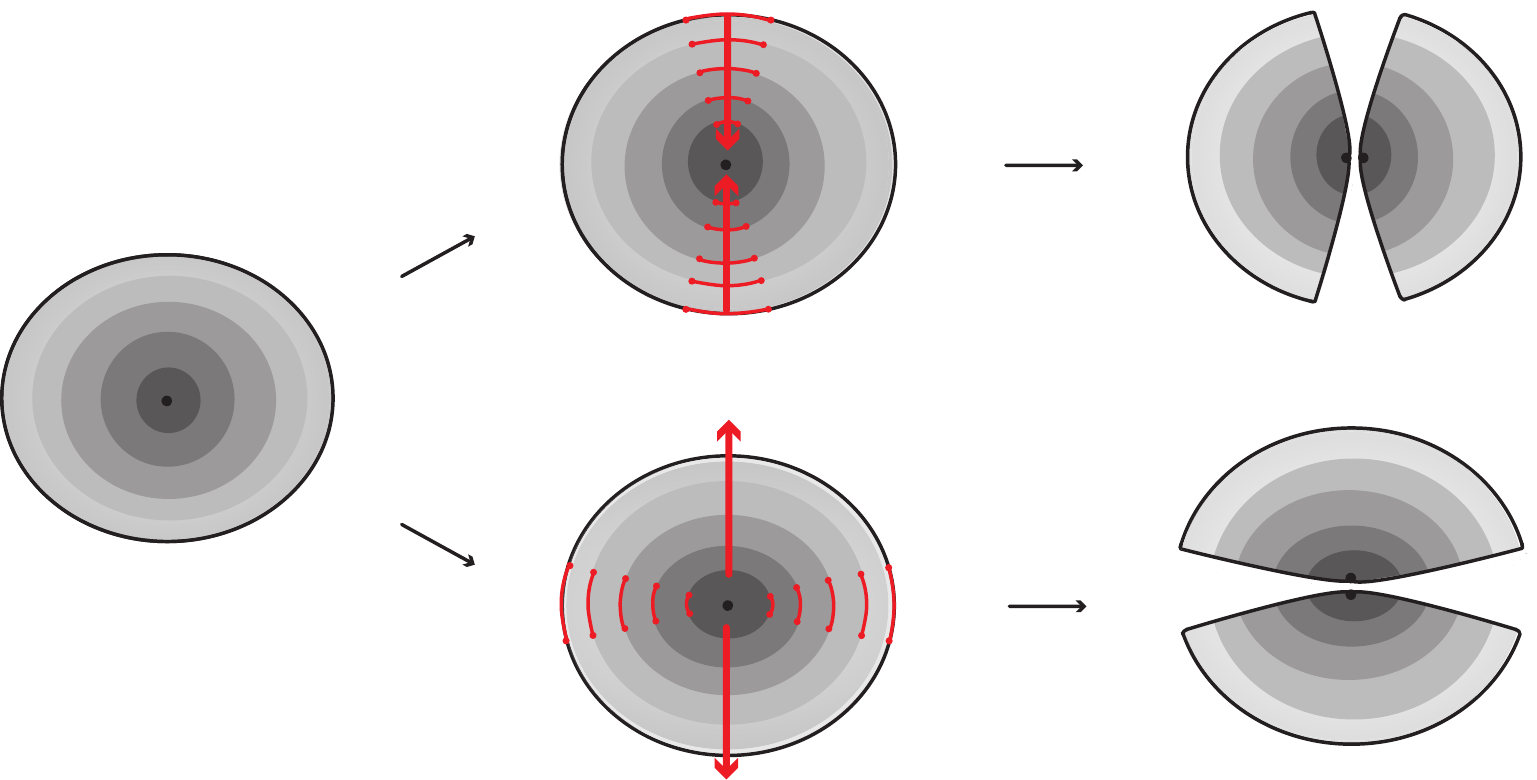}
\caption{Solid 1-dimensional surgery.}
\label{solid1D}
\end{center}
\end{figure}

There are phenomena which seem to undergo the process of  1-dimensional surgery but happen on surfaces, such as tension on membranes or soap films. In order to model topologically such phenomena we introduce the notion of solid 1-dimensional surgery. {\it Solid 1-dimensional surgery} on the $2$-disc, $D^2$, is the topological procedure whereby a ribbon $D^1 \times D^1$ is being removed, such that the closure of the remaining manifold comprises two discs  $D^2 \times S^0$. See Figure~\ref{Formal1D} where the interior is now supposed to be filled in.
 This process is equivalent to performing 1-dimensional surgeries on the whole continuum of concentric circles included in $D^2$. More precisely, and introducing at the same time dynamics, we define:

\begin{defn} \rm  We start with the $2$-disc of radius 1 with polar layering: 
$$
D^2 = \cup_{0<r\leq 1} S^1_r \cup \{C\},
$$ 
where $r$ the radius of a circle and $C$ the limit point of the circles, that is, the center of the disc. We specify colinear pairs of antipodal points, with neighbourhoods of analogous lengths, on which the same colinear forces act, attracting or repelling, see Figure~\ref{solid1D}. Then we perform  1-dimensional surgery, attracting or repelling, on the whole continuum of  concentric circles. We also  define 1-dimensional surgery on the limit point $C$ to be the two limit points of the resulting surgeries. That is, the effect of {\it  1-dimensional surgery on a point is the creation of two new points}. The above process is the same as first removing the center $C$ from $D^2$, doing the  1-dimensional surgeries and then taking the closure of the resulting space, see Figure~\ref{solid1D}. The resulting manifold is 
$$
\chi(D^2) := \cup_{0<r\leq 1}\chi(S^1_r) \cup \chi(c),
$$
which comprises two copies of $D^2$. {\it Attracting solid 1-dimensional surgery} on $D^2$ is the above topological procedure whereby attracting forces act on the circles $ S^1_r$, see Figure~\ref{solid1D}. {\it Repelling solid 1-dimensional surgery} on $D^2$ is the above topological procedure whereby repelling forces act on the circles $ S^1_r$, see Figure~\ref{solid1D}.  
\end{defn}

\smallbreak

\section{2-dimensional topological surgery}\label{2D}

\subsection{}

\begin{figure}[ht!]
\begin{center}
\includegraphics[width=14.5cm]{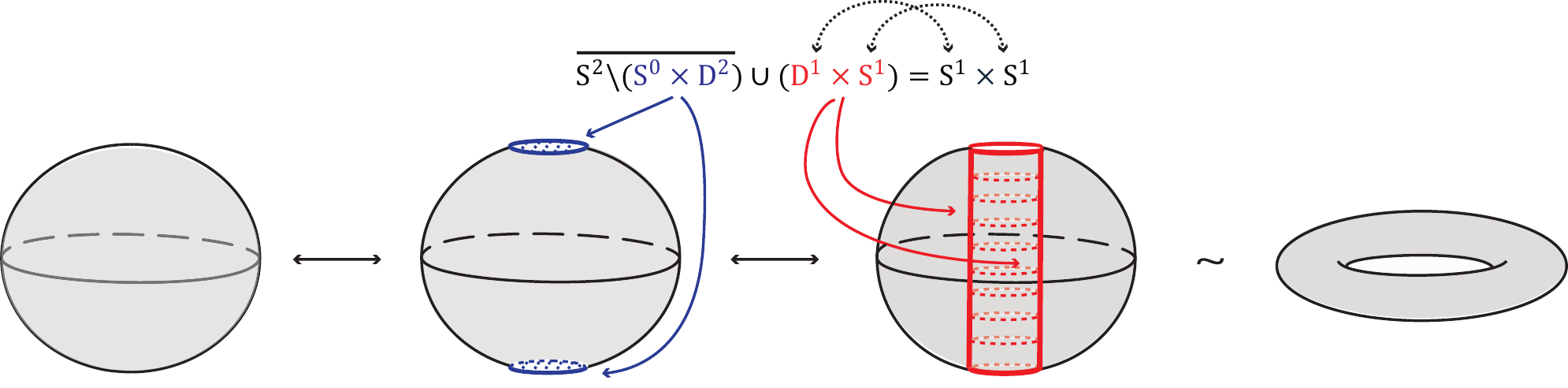}
\caption{Surgery on the sphere results in a torus.}
\label{formal2Da}
\end{center}
\end{figure}

Starting with $S^2$, there are two types of {\it 2-dimensional topological surgery}. In the first type two discs $S^0\times D^2$ are removed from  $S^2$ and are replaced in the closure of the remaining manifold by a cylinder $D^1\times S^1$, which gets attached along the common boundary $S^0\times S^1$, comprising two copies of $S^1$, via a homeomorphism. The gluing homeomorphism of the common boundary is a number of full twists of each copy of $S^1$. The above operation changes the homeomorphism type from the 2-sphere to that of the torus (see Figure~\ref{formal2Da}). In fact, every c.c.o. surface  arises from the 2-sphere by repeated surgeries and each time the above process is performed the genus of the surface is increased by one.  Note that, if the cylinder were attached on $S^2$ externally, the result would still be a torus. Physical examples reminiscent of 2-dimensional surgery comprise the formation of whirls and the Falaco solitons \cite{Ki} (see Figure~\ref{falaco}).

\smallbreak

\begin{figure}[ht!]
\begin{center}
\includegraphics[width=14cm]{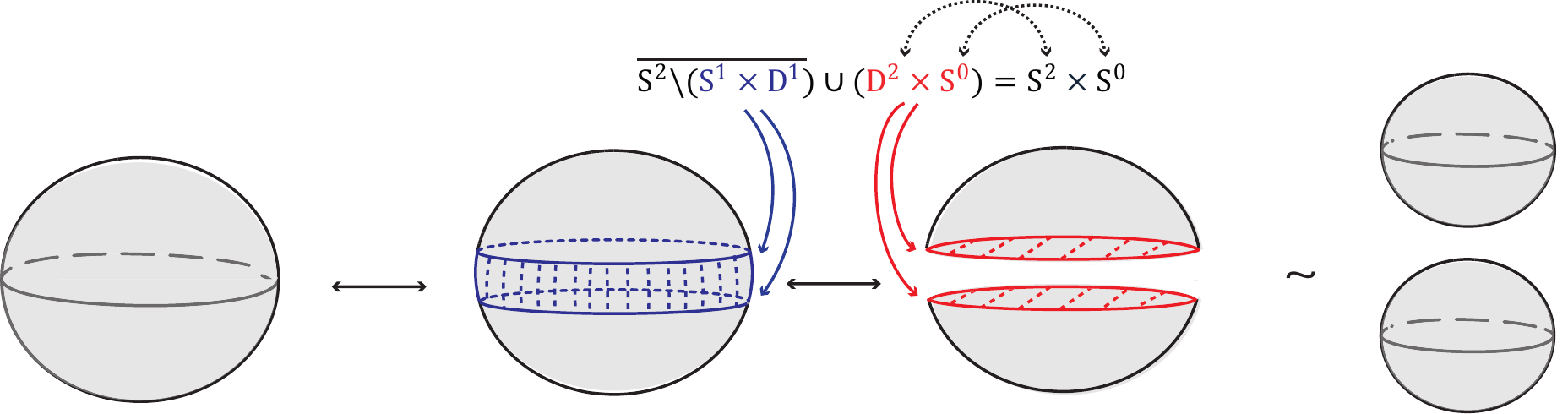}
\caption{Surgery on the sphere results in two spheres.}
\label{formal2Dr}
\end{center}
\end{figure}

\begin{figure}[ht!]
\begin{center}
\includegraphics[width=10cm]{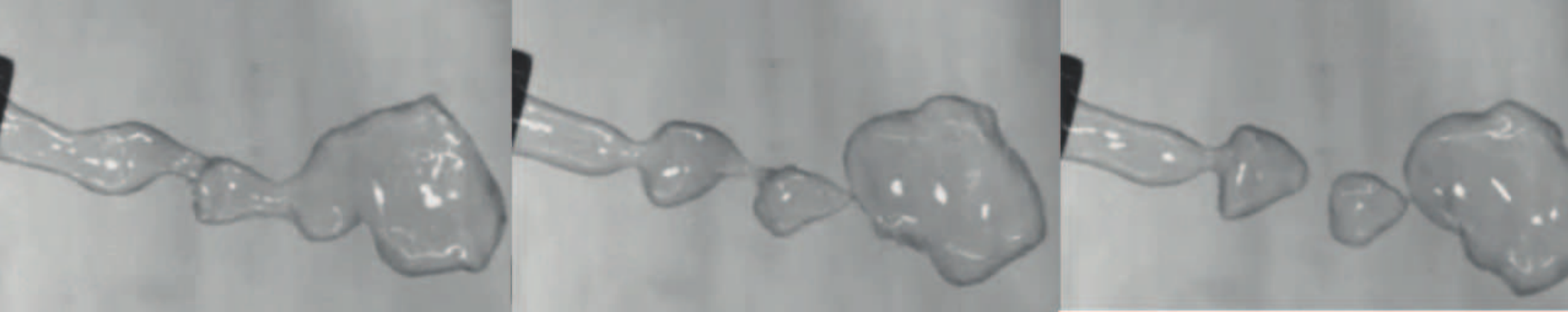}
\caption{Soap bubble blowing.}
\label{Bubbles}
\end{center}
\end{figure}

The other possibility of 2-dimensional surgery on the 2-sphere is the following:  an annulus $S^1 \times D^1$ (perhaps twisted a number of times) is removed from $S^2$ and is replaced in the closure of the remaining manifold by two discs  $ D^2 \times S^0$ attached along the common boundary $S^1 \times S^0$, resulting in two copies of  $S^2$. See Figure~\ref{formal2Dr}. Phenomena exemplifying this type of surgery comprise soap bubble blowing and, similarly, glass blowing, see Figure~\ref{Bubbles}.  
 It is worth noting that this type of surgery applied on a torus is the reverse process of the attracting type. Namely, if  a cylinder were removed from a torus and were replaced by two discs the result would be a 2-sphere. 

\smallbreak

\subsection{}

In order to model topologically phenomena exhibiting 2-dimensional surgery or to understand 2-dimensional surgery through continuity we need, also here, to introduce dynamics. 

\begin{figure}[ht!]
\begin{center}
\includegraphics[width=15cm]{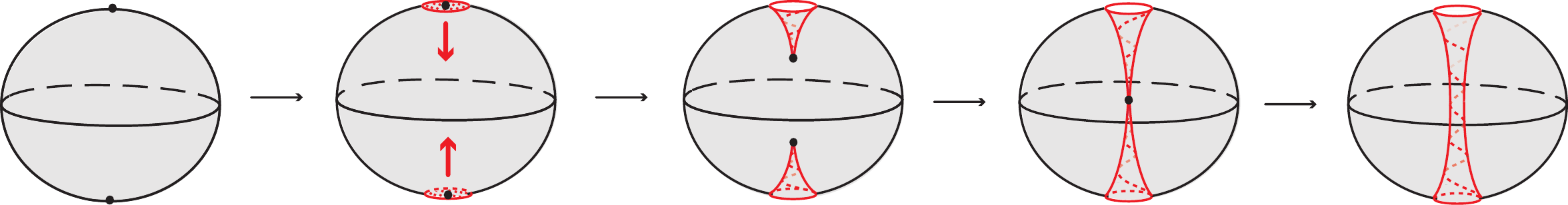}
\caption{Attracting 2-dimensional surgery.}
\label{2DAttract}
\end{center}
\end{figure}

\begin{figure}[ht!]
\begin{center}
\includegraphics[width=14.5cm]{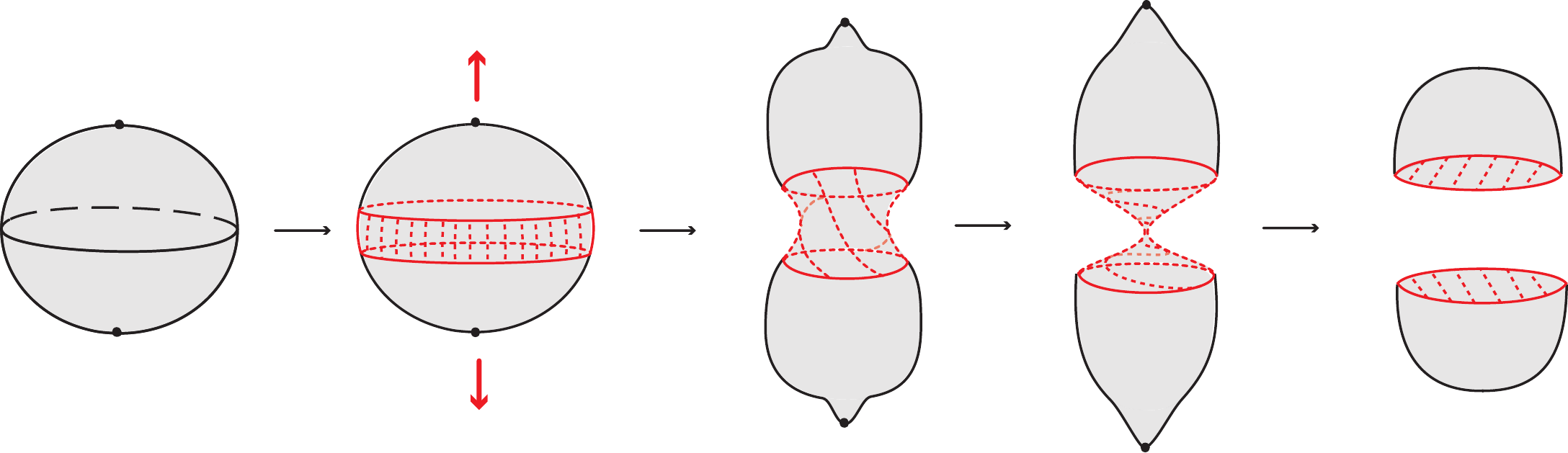}
\caption{Repelling 2-dimensional surgery.}
\label{2DRepel}
\end{center}
\end{figure}

\begin{defn} \rm The {\it attracting 2-dimensional surgery} starts with two  poles specified on $S^2$ with attracting forces applied on them.  Then two discs $S^0\times D^2$,  neighbourhoods of the two poles, approach each other, with a possible number of full twists. When the two discs touch, {\it recoupling} takes place and the discs get transformed into the final cylinder. See Figure~\ref{2DAttract}. The twisting propagates along the cylinder, reminding the process of hole drilling.

 In the {\it repelling 2-dimensional surgery} two  poles are specified on $S^2$ with repelling forces pulling them to opposite directions. This creates, by passive reaction, a cylindrical `necking' in the middle, which eventually tears apart and new material,  two discs, gets attached along the boundary $S^1 \times S^0$. See Figure~\ref{2DRepel}. 
\end{defn}

\smallbreak

\smallbreak

\smallbreak

\begin{rem} \rm It is worth observing that the  process of repelling 2-dimensional surgery in reverse time would mean that the initial surface comprises two copies of $S^2$ and there are two discs to be removed, one on each sphere, replaced by a cylinder, thus merging the two spheres into one. Similarly,  the  process of attracting 2-dimensional surgery in reverse time would mean that the initial surface is 
the torus and there is a cylinder to be removed and replaced by two discs, thus yielding back the 2-sphere. In other words, the reverse process of repelling surgery (where repelling forces are applied on the boundary circles) can be viewed as attracting surgery (where the attracting forces are now applied on the centers of the two discs) and vice versa. 
\end{rem}

\subsection{} 

In some natural phenomena the object undergoing surgery is not a surface but three-dimensional. For this reason we introduce also here the notion of {\it solid 2-dimensional surgery}. There are two types of solid 2-dimensional surgery on the $3$-ball, $D^3$, analogous to the two types of 2-dimensional surgery.
The first one is the topological procedure of removing a solid cylinder homeomorphic to the product set $D^1 \times D^2$, $h(D^1 \times D^2)$ (such that the part $S^0 \times D^2$ of its boundary lies in the boundary of $D^3$) and taking the closure of the remaining manifold $D^3 \setminus h(D^1 \times D^2)$, which is a twisted solid torus.  See Figure~\ref{formal2Da} where the interior is supposed to be filled in.
 The second type is the topological procedure of removing a solid cylinder homeomorphic to the product set $D^2 \times D^1$, $h(D^2 \times D^1)$, (such that the part $S^1 \times D^1$ of its boundary lies in the boundary of $D^3$) and taking the closure of the remaining manifold $D^3 \setminus h(D^2 \times D^1)$, which is two copies of $D^3$.  See Figure~\ref{formal2Dr} where the interior is supposed to be filled in. 
 
\subsection{}

\begin{figure}[ht!]
\begin{center}
\includegraphics[width=15cm]{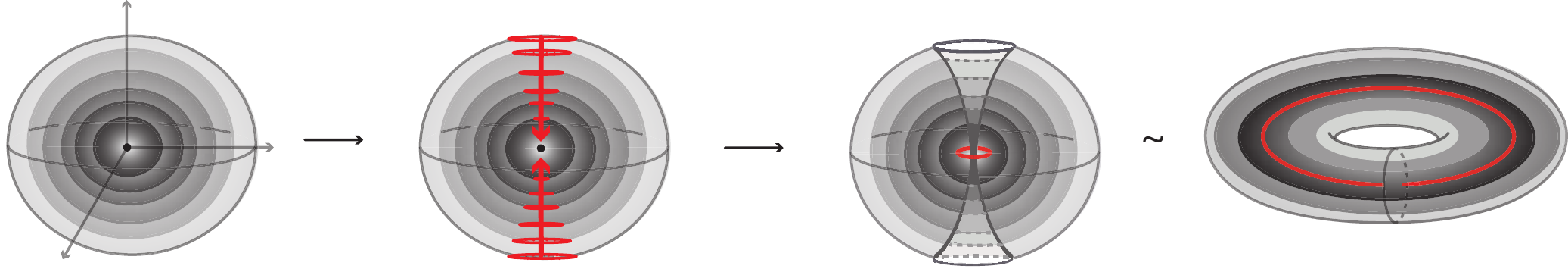}
\caption{Attracting solid 2-dimensional surgery.}
\label{solid2DAttract}
\end{center}
\end{figure}

\begin{figure}[ht!]
\begin{center}
\includegraphics[width=11cm]{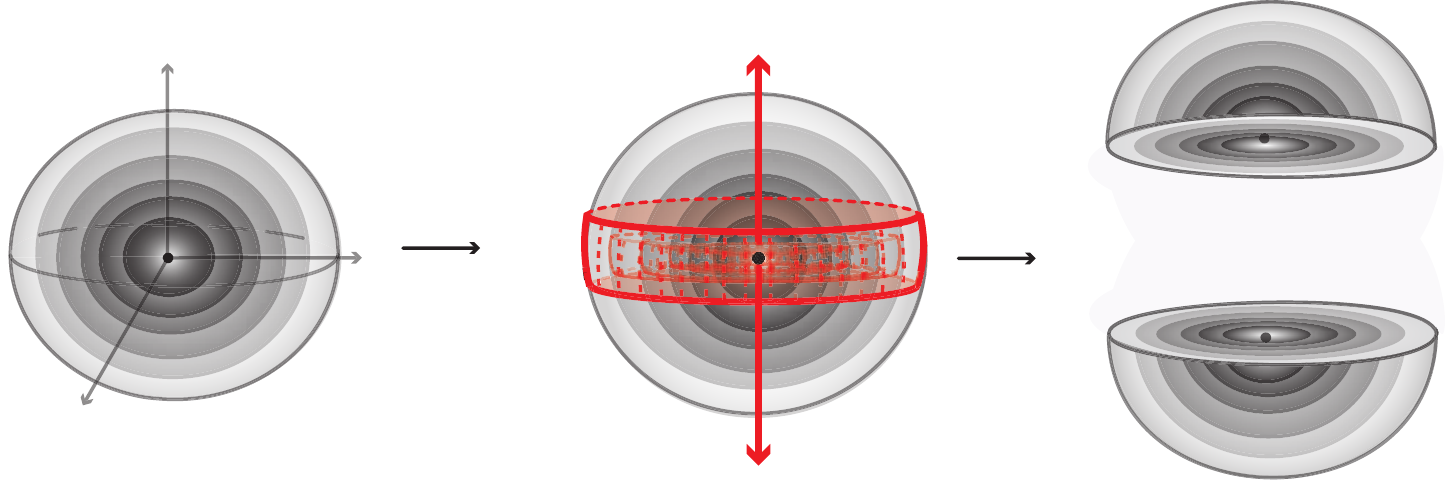}
\caption{Repelling solid 2-dimensional surgery.}
\label{solid2DRepel}
\end{center}
\end{figure}

In order to model better natural phenomena exemplifying solid 2-dimensional surgery we  shall introduce  dynamics:

\begin{defn}\label{continuum2D} \rm  Start with the $3$-ball of radius 1 with polar layering: 
$$
D^3 = \cup_{0<r\leq 1} S^2_r \cup \{C\},
$$ 
where $r$ the radius of a 2-sphere and $C$ the limit point of the spheres, that is, the center of the ball. 
{\it Attracting solid 2-dimensional surgery} on $D^3$ is the  topological procedure where: on all spheres $S^2_r$  colinear pairs of antipodal points are specified, on which the same colinear  attracting forces act. The poles have disc neighbourhoods of analogous areas, see Figure~\ref{solid2DAttract}. Then  attracting 2-dimensional surgeries are performed on the whole continuum of the concentric spheres using the same homeomorphism $h$. Moreover,  attracting 2-dimensional surgery on the limit point $C$ is defined to be the  limit circle of the nested tori resulting from the continuum of 2-dimensional surgeries, see Figure~\ref{solid2DAttract}. That is, the effect {\it  of attracting 2-dimensional surgery on a point is the creation of a circle}.  The process is characterized by the 1-dimensional core $L$ of the removed solid cylinder, joining the antipodal points on the outer shell and intersecting each spherical layer in the two antipodal points; also, by the homeomorphism $h$, resulting in the whole continuum of layered tori, and it can be viewed as drilling out a tunnel along $L$ according to $h$.  For $h$ non-trivial, this agrees with our intuition that, for opening a hole, drilling with twisting seems to be the easiest way.
\smallbreak
 {\it Repelling solid 2-dimensional surgery} on $D^3$  is the  topological procedure where: on all spheres $S^2_r$ nested  annular peels of the solid annulus of analogous areas are specified and the same colinear repelling forces act on all spheres, see Figure~\ref{solid2DRepel}. 
Then  repelling 2-dimensional surgeries are performed on the whole continuum of the concentric spheres using the same homeomorphism $h$, see Figure~\ref{solid2DRepel}. Moreover, repelling 2-dimensional surgery  on the limit point $C$  is defined to be the two limit points of the nested pairs of 2-spheres resulting from  the continuum of 2-dimensional surgeries, see Figure~\ref{solid2DRepel}. That is, the effect of {\it  repelling 2-dimensional surgery on a point is the creation of two new points}. 
The process is characterized by the 2-dimensional central disc of the solid annulus and the homeomorphism $h$, and it can be viewed as pulling apart along the central disc, after a number of twists according to  $h$. For $h$ non-trivial, this operation agrees with our intuition that for cutting a solid object apart, pulling with twisting seems to be the easiest way.

 In either case the above process is the same as first removing the center $C$ from $D^3$, performing the  2-dimensional surgeries and then taking the closure of the resulting  space. Namely we obtain:  
$$
\chi(D^3) := \cup_{0<r\leq 1}\chi(S^2_r) \cup \chi(C),
$$
which is a solid torus in the case of attracting solid 2-dimensional surgery and  two copies of $D^3$ in the case of repelling solid 2-dimensional surgery. See Figures~\ref{solid2DAttract} and~\ref{solid2DRepel}.
\end{defn}

\smallbreak

\smallbreak

\begin{note} \rm The notions of  2-dimensional (resp. solid 2-dimensional) surgery, attracting or repelling, can be generalized from $S^2$ (resp. $D^3$) to a surface (resp. handlebody of genus $g$) creating a surface (resp. handlebody of genus $g+1$).
\end{note} 

\begin{figure}[ht!]
\begin{center}
\includegraphics[width=15cm]{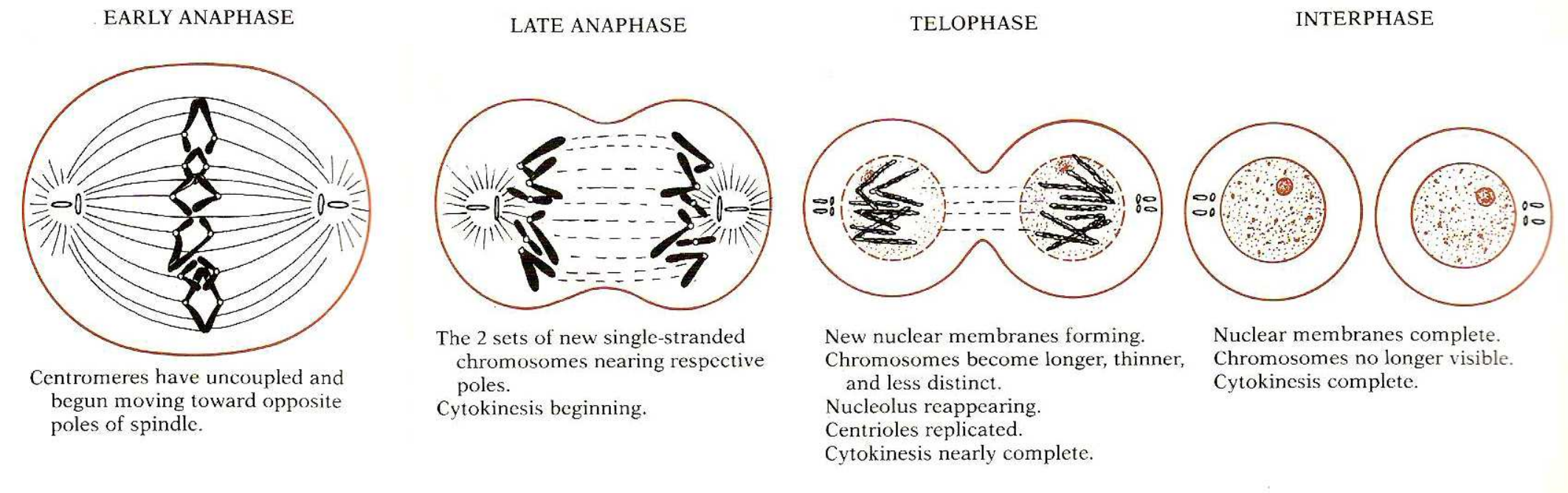}
\caption{The process of mitosis is an example of solid repelling 2-dimensional surgery.}
{\footnotesize \underline{{\it Source}}: W.T. Keeton, C.H. McFadden, {\em Elements of Biological Science}, W.W. Norton \& Company Inc., 3rd edition (1983), p. 395.}
\label{mitosis}
\end{center}
\end{figure}

\begin{figure}[ht!]
\begin{center}
\includegraphics[width=4cm]{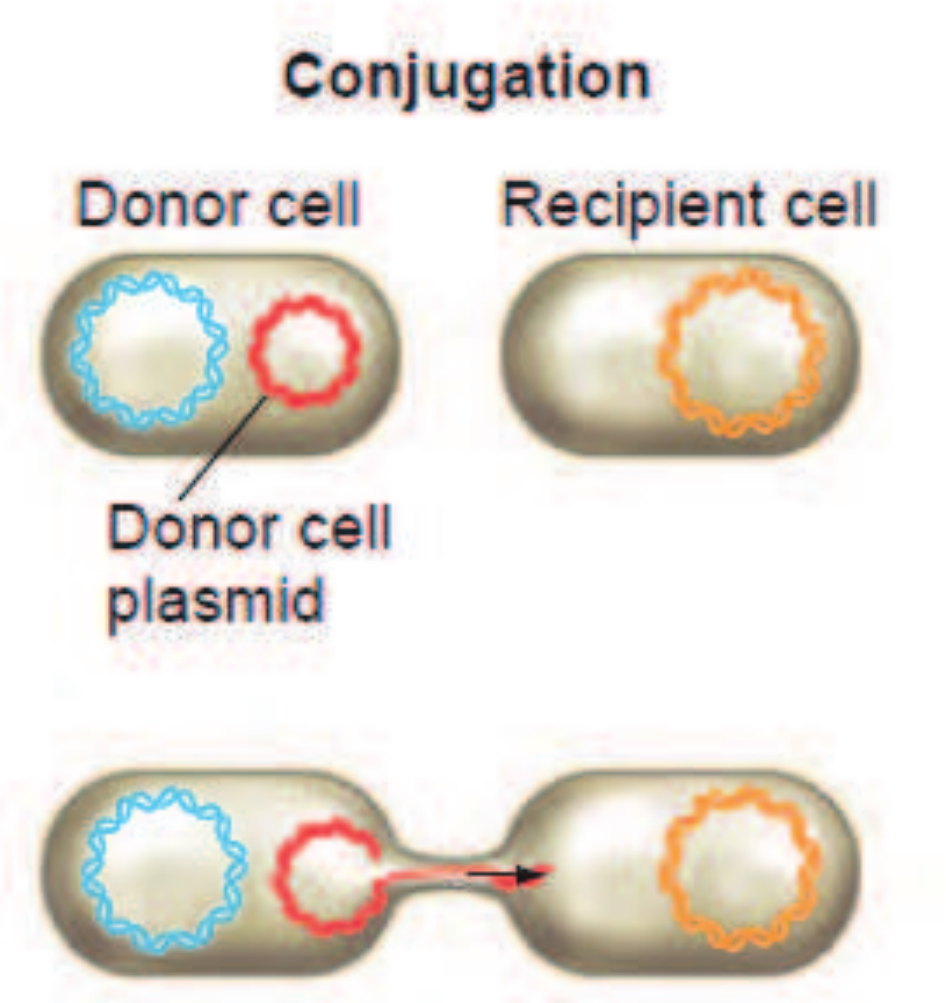}
\caption{Gene transfer in bacteria.}
{\footnotesize \underline{{\it Source}}: \cite{HHGRSV}, p. 486.}
\label{genetransfer}
\end{center}
\end{figure}

A good natural model reminiscent of attracting solid 2-dimensional surgery  is the formation of an apple (formed from a seed). Here no twisting occurs, so $h$ is trivial. Repelling solid 2-dimensional surgery can be exemplified by the biological process of mitosis, where a cell splits into two new cells. See Figure~\ref{mitosis} (for description and instructive illustrations see for example \cite{KeFa}, p. 395). Further, it is worth noting that the reverse process of  repelling solid 2-dimensional surgery can be  found in the mechanism of gene transfer in bacteria. See Figure~\ref{genetransfer}  (for description and instructive illustrations see, for example, \cite{HHGRSV}). Here ``donor DNA is transferred directly to recipient through a connecting tube" and two copies of $D^3$ merge in one.

\smallbreak

\smallbreak

\subsection{}

\begin{figure}[ht!]
\begin{center}
\includegraphics[width=10cm]{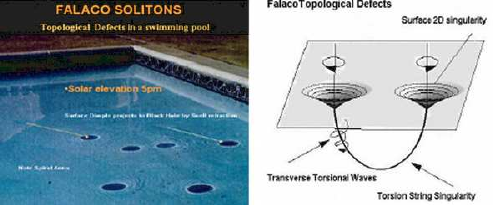}
\caption{Three pairs of Falaco Solitons in a swimming pool.}
{\footnotesize \underline{{\it Source}}: R.M. Kiehn, {\em Non-Equilibrium Systems and Irreversible Processes}, Adventures in Applied  Topology {\bf 1}, Non Equilibrium Thermodynamics, University of Houston Copyright CSDC. INC, (2004), pp. 147, 150.}
\label{falaco}
\end{center}
\end{figure}

Attracting solid 2-dimensional surgery can be also observed in the formation of the Falaco solitons \cite{Ki} (Figure~\ref{falaco}) and in the formation of whirls.  The Falaco solitons are pairs of singular surfaces (poles) connected by means of a stabilizing invisible thread. As shown in Figure~\ref{falaco}, starting by the two poles, and by drilling along the joining line, surgery seems to be performed. Based on the experimental creation of Falaco solitons in a swimming pool, it has been conjectured that M31 and the Milky Way galaxies could be connected by a `topological thread'.

\smallbreak

In such phenomena we do not see the whole space $D^3$; the process can be viewed as taking place between the two attracting discs of the solid cylinder, so the initial space can be considered to be $D^1 \times D^2 = \overline{D^3 \setminus D^2 \times D^1}$. This type of surgery shall be called {\it truncated attracting solid 2-dimensional surgery}.  In the above examples $h$ is non-trivial. 

\smallbreak

One could also define theoretically the non-solid analogue, the {\it truncated attracting 2-dimensional surgery} as attracting 2-dimensional surgery taking place just between the two attracting discs, which are neighbourhoods of the two specified points on $S^2$. So, the initial manifold can be considered to be just these two discs, that is, $S^0 \times D^2 = \overline{S^2 \setminus S^1 \times D^1}$.

\begin{figure}[ht!]
\begin{center}
\includegraphics[width=4.7cm]{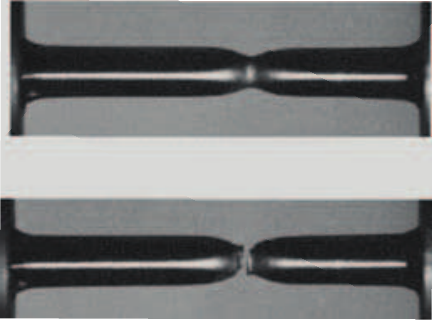}
\caption{Tension and the necking phenomenon.}
{\footnotesize \underline{{\it Source}}: http://www.ara.com/Projects/SVO/weld.htm.}
\label{necking}
\end{center}
\end{figure}

\smallbreak
Another phenomenon falling topologically in the scheme of repelling solid 2-dimensional surgery is tension on metal speciments and the `necking effect'. More precisely, in experiments in mechanics tensile forces (or loading)  are applied on a cylindrical speciment made of dactyle material (steel, aluminium, etc.). Up to some critical value of the force the deformation is homogenuous (the cross-sections have the same area). At the critical value the deformation is localized within a very small area where the cross-section is reduced drastically, while the sections of the remaining portions increase slightly. This is the `necking phenomenon'. Shortly after the speciment is fractured. View Figure~\ref{necking}.

\smallbreak

In such phenomena we do not see the whole space $D^3$; the process can be seen as being localized just in the region of the solid annulus, so the initial space can be considered to be $D^2 \times D^1$. This type of surgery shall be called  {\it truncated repelling solid 2-dimensional surgery}. One could also define theoretically the non-solid analogue, the {\it truncated repelling 2-dimensional surgery} as repelling 2-dimensional surgery taking place just in the region of the annulus $S^1 \times D^1$ which is complementary to the two repelling discs. So, the initial manifold can be considered to be just this annulus, that is, $S^1 \times D^1 = S^2 \setminus S^0 \times D^2$.

\begin{rem} \rm A cross-section of 2-dimensional surgery of attracting  or repelling type, truncated or solid, passing through the specified points is precisely the corresponding type of 1-dimensional surgery.
\end{rem}

\section{3-dimensional topological surgery}\label{3D}

\subsection{}

In dimension 3, the simplest c.c.o. 3-manifolds are: the 3-sphere $S^3$ and the lens spaces $L(p,q)$.  We start with  $S^3$ and we recall its three most common descriptions. 

\smallbreak

\begin{figure}[ht!]
\begin{center}
\includegraphics[width=12cm]{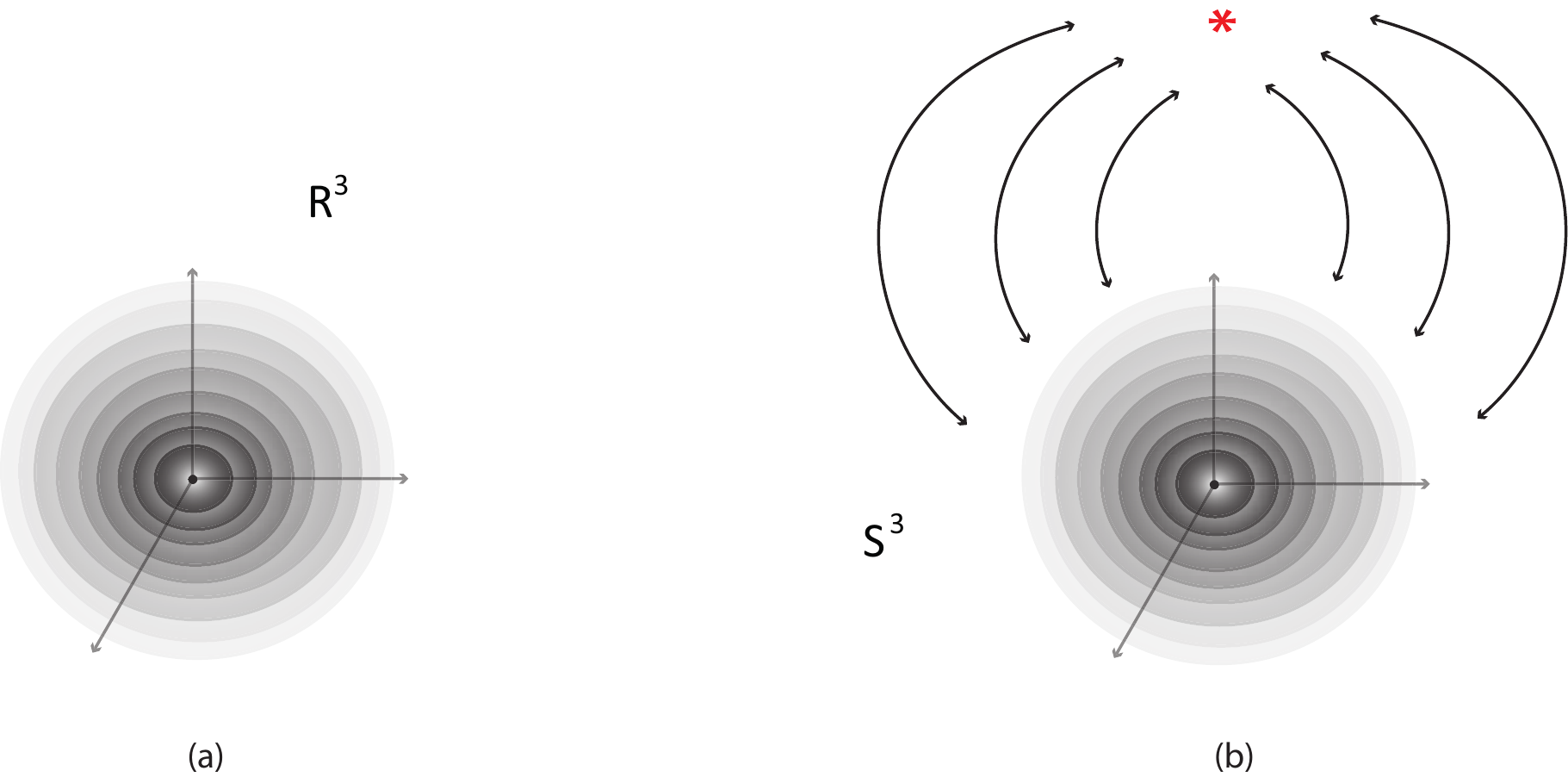}
\caption{$S^3$ is the compactification of ${\mathbb R}^3$.}
\label{layeredspheres}
\end{center}
\end{figure}

Firstly, $S^3$ can be viewed as ${\mathbb R}^3$ with all points at infinity compactified to one single point: $S^3 = {\mathbb R}^3 \cup \{\infty\}$. See Figure~\ref{layeredspheres}(b). 
 ${\mathbb R}^3$ can be viewed as an unbounded continuum of nested 2-spheres centered at the origin, together with the point at the origin, see Figure~\ref{layeredspheres}(a),  and also as the de-compactification of $S^3$.  So, $S^3$ minus the point at the origin and the point at infinity can be viewed as a continuous nesting of 2-spheres.   

\smallbreak

\begin{figure}[ht!]
\begin{center}
\includegraphics[width=12cm]{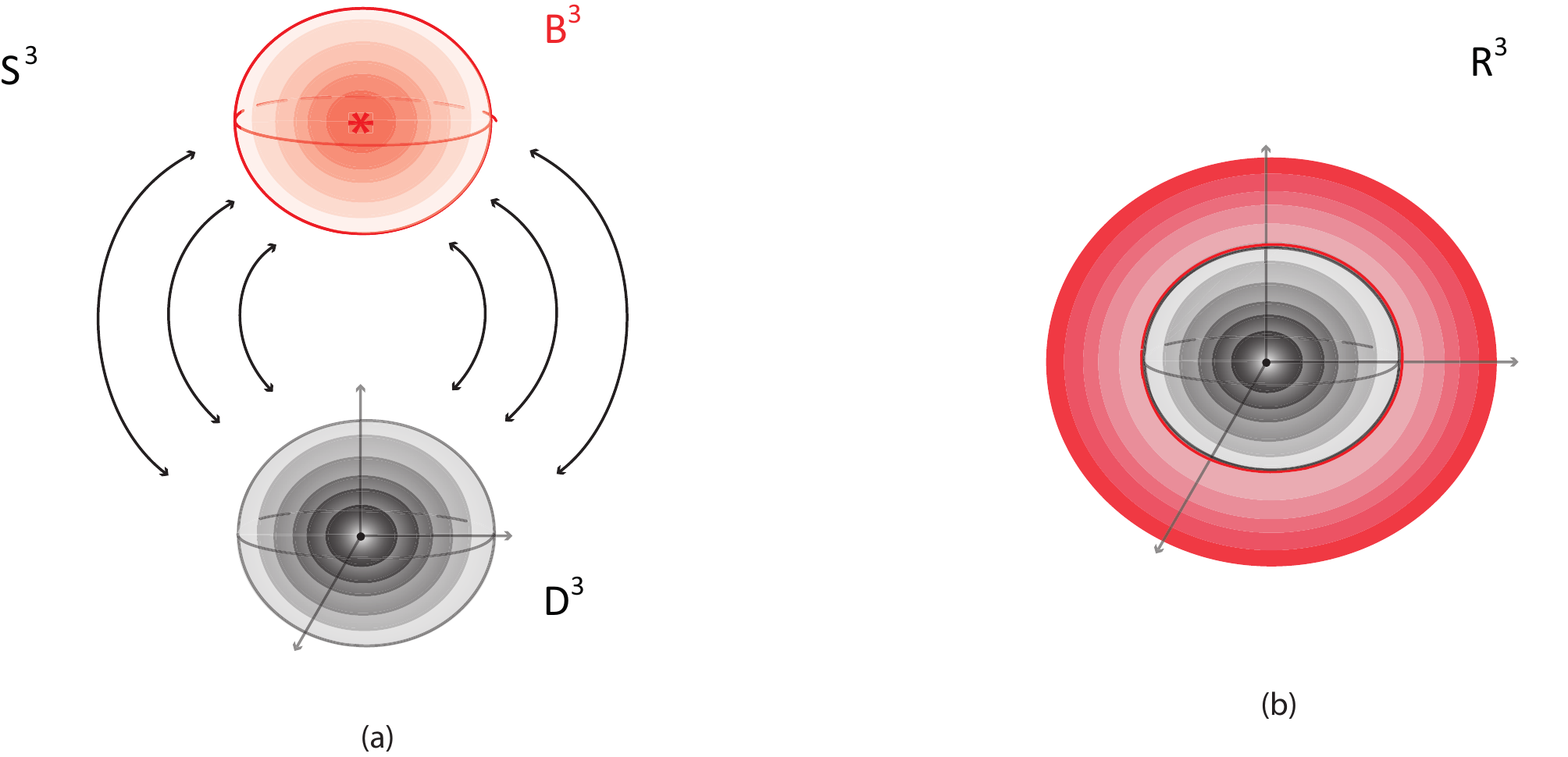}
\caption{$S^3$ is the result of gluing two 3-balls.}
\label{twoballs}
\end{center}
\end{figure}

Secondly, $S^3$ can be viewed as  the union of two $3$-balls: $S^3 = B^3 \cup D^3$, see Figure~\ref{twoballs}(a). The two descriptions of $S^3$ are clearly related, since a  (closed) neighbourhood of the point at infinity can stand for one of the two $3$-balls.  Note that, when removing the point at infinity in Figure~\ref{twoballs}(a) we can see the concentric spheres of the 3-ball $B^3$ (in red) wrapping around the concentric spheres of the 3-ball $D^3$, see Figure~\ref{twoballs}(b). This is another way of viewing ${\mathbb R}^3$ as the de-compactification of $S^3$.  This picture is the analogue of the stereographic projection of $S^2$ on the plane ${\mathbb R}^2$, whereby the projections of the concentric circles of the south hemisphere together with the projections of the concentric circles of the north hemisphere form the well-known polar description of ${\mathbb R}^2$ with the unbounded continuum of concentric circles. 

\smallbreak

\begin{figure}[ht!]
\begin{center}
\includegraphics[width=11cm]{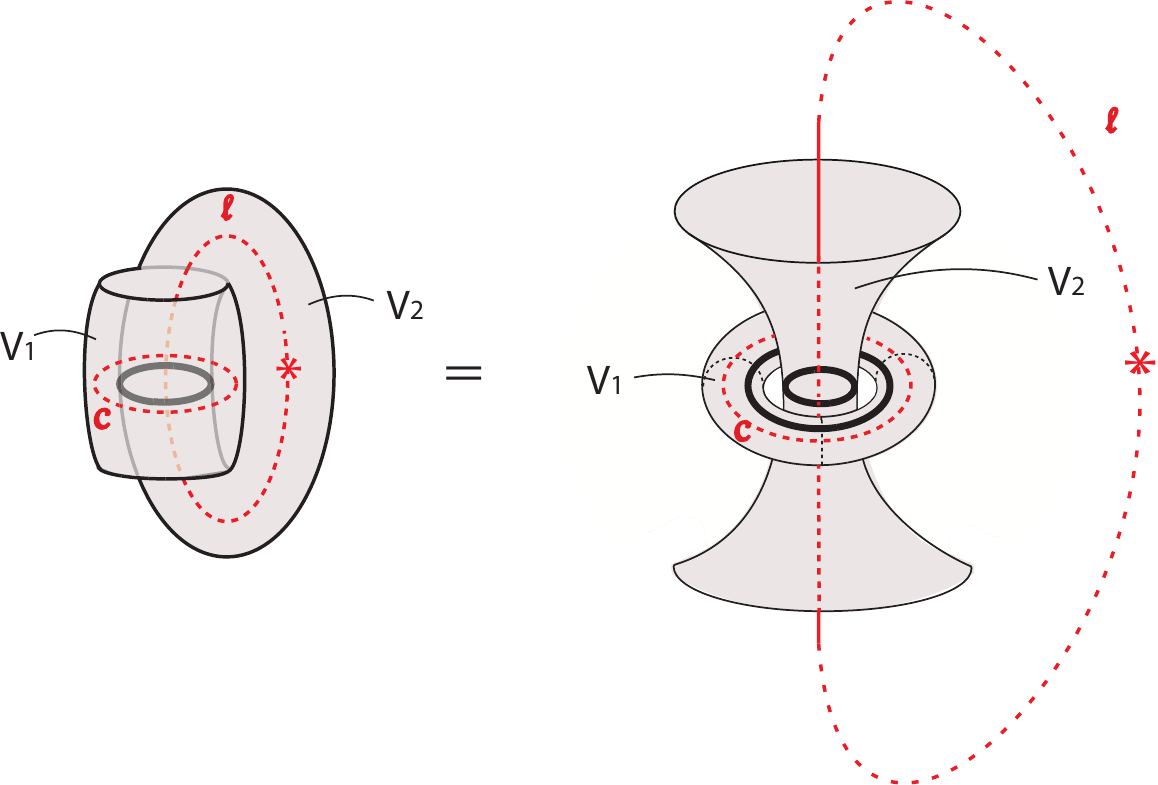}
\caption{$S^3$ as a union of two solid tori.}
\label{SplittingofS3}
\end{center}
\end{figure}

The third well-known representation of $S^3$ is as the union of two solid tori,  $S^3 = V_1\cup_\vartheta V_2$, via the torus homeomorphism $\vartheta$ along the common boundary.  $\vartheta$ maps a meridian of $V_2$ to a longitude of $V_1$ which has linking number zero with the core curve $c$ of $V_1$.
The illustration in Figure~\ref{SplittingofS3} gives an idea of this splitting of $S^3$. In the figure, the core curve of $V_1$ is in dashed red. So, the complement of a solid torus $V_1$ in $S^3$ is another solid torus $V_2$ whose core curve $l$ (the dashed red curve in the figure)  may be assumed to pass by the point at infinity. Note that, $S^3$ minus the core curves $c$ and $l$ of  $V_1$ and $V_2$ (the red curves in Figure~\ref{SplittingofS3}) can be  viewed as a continuum of nested tori.   

\smallbreak

\begin{figure}[ht!]
\begin{center}
\includegraphics[width=9cm]{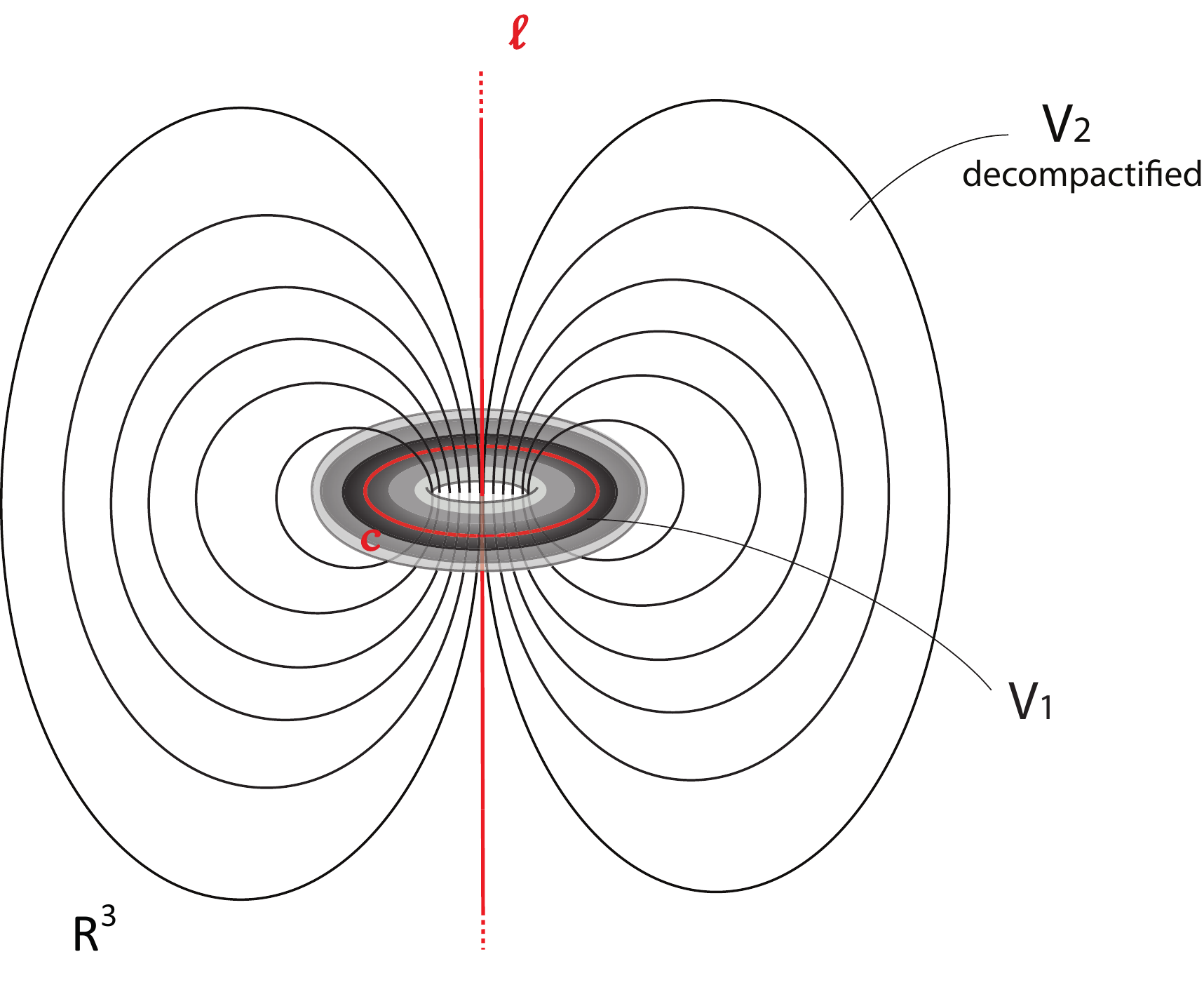}
\caption{De-compactification of $S^3$ viewed as two tori.}
\label{DecompS3}
\end{center}
\end{figure}

When removing the point at infinity in the representation of $S^3$ as a union of two solid tori, the core of the solid torus $V_2$ becomes an infinite line $l$ and the nested tori of  $V_2$ can now be seen wrapping around the nested tori of  $V_1$. See Figure~\ref{DecompS3}.  Therefore, ${\mathbb R}^3$ can be viewed as an unbounded continuum of nested tori, together with the core curve $c$ of  $V_1$  and the infinite line  $l$. This line $l$ joins pairs of antipodal points of all concentric spheres of the first description. Note that in the nested spheres description (Figure~\ref{layeredspheres}) the  line $l$ pierces all spheres while in the nested tori description the  line $l$ is the `untouched' limit circle of all tori.

\smallbreak

\subsection{}

\begin{figure}[ht!]
\begin{center}
\includegraphics[width=13cm]{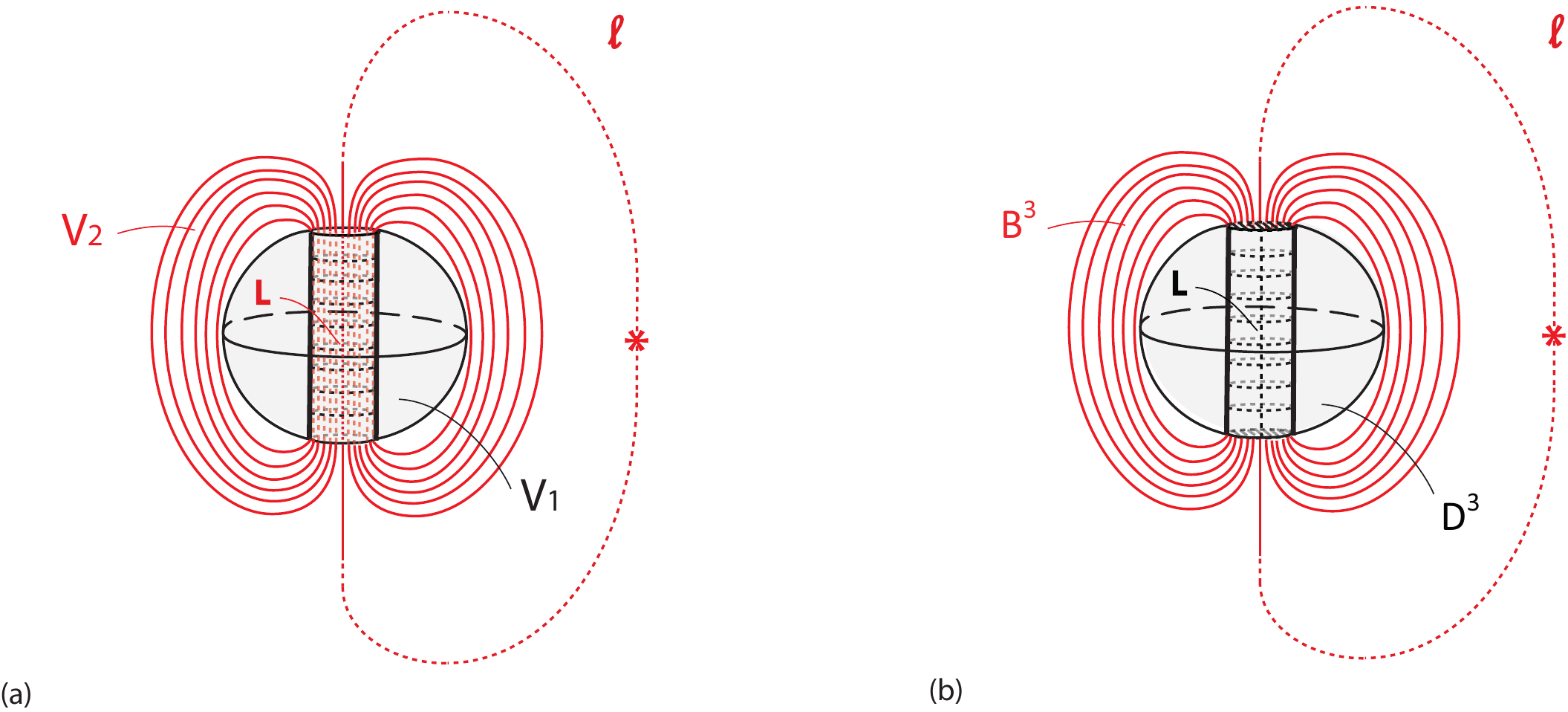}
\caption{Passing from (a) $S^3$ as two tori  to (b) $S^3$ as two balls.}
\label{S3BallsToTori}
\end{center}
\end{figure}

The third  description of $S^3$ is a bit harder to connect with the first two. We shall do this here. A way to see this connection is  the following. 
Consider the description of  $S^3$ as the union of two 3-balls, $B^3$ and $D^3$ (Figure~\ref{layeredspheres}(b)). 
Combining with the third description of  $S^3$ (Figure~\ref{SplittingofS3}) we notice that both 3-balls are pierced by  the core curve $l$ of the solid torus $V_2$. 
 Therefore, $D^3$ can be viewed as the solid torus $V_1$ to which  a solid cylinder $D^1\times D^2$  is attached  via the homeomorphism $\vartheta$: 
$$
D^3 = V_1\cup_\vartheta (D^1\times D^2).
$$
This solid cylinder is part of the solid torus $V_2$,  a `cork'  filling the hole of $V_1$. Its core curve is an arc $L$, part of the core curve $l$ of $V_2$. View  Figure~\ref{S3BallsToTori}. The second ball $B^3$ (Figure~\ref{layeredspheres}(b)) can be viewed as the remaining of $V_2$ after removing the cork $D^1\times D^2$: 
$$
B^3 = \overline{ V_2 \setminus_\vartheta (D^1\times D^2)}.
$$
In other words the solid torus $V_2$ is cut into two solid cylinders, one comprising the `cork' of $V_1$ and the other comprising the 3-ball $B^3$. 

\smallbreak

\begin{rem}\label{truncate} \rm
If we remove a whole neighbourhood $B^3$ of the point at infinity and focus on the remaining  3-ball $D^3$, the line $l$ of the previous picture is truncated to the arc $L$ and the solid cylinder $V_2$ is truncated to  the cork of $D^3$.   
\end{rem}

Another way to see the connection among the different descriptions of  $S^3$ is by 
combining the above with Definition~\ref{continuum2D}. Indeed, one can pass from the second description of $S^3$ to the third by performing attracting solid 2-dimensional surgery (with trivial homenomorphism)  on the 3-ball $D^3$ along the arc $L$. Note that, by  Definition~\ref{continuum2D}, the point at the origin turns into the core curve of $V_1$.

\subsection{}

\begin{figure}[ht!]
\begin{center}
\includegraphics[width=8.5cm]{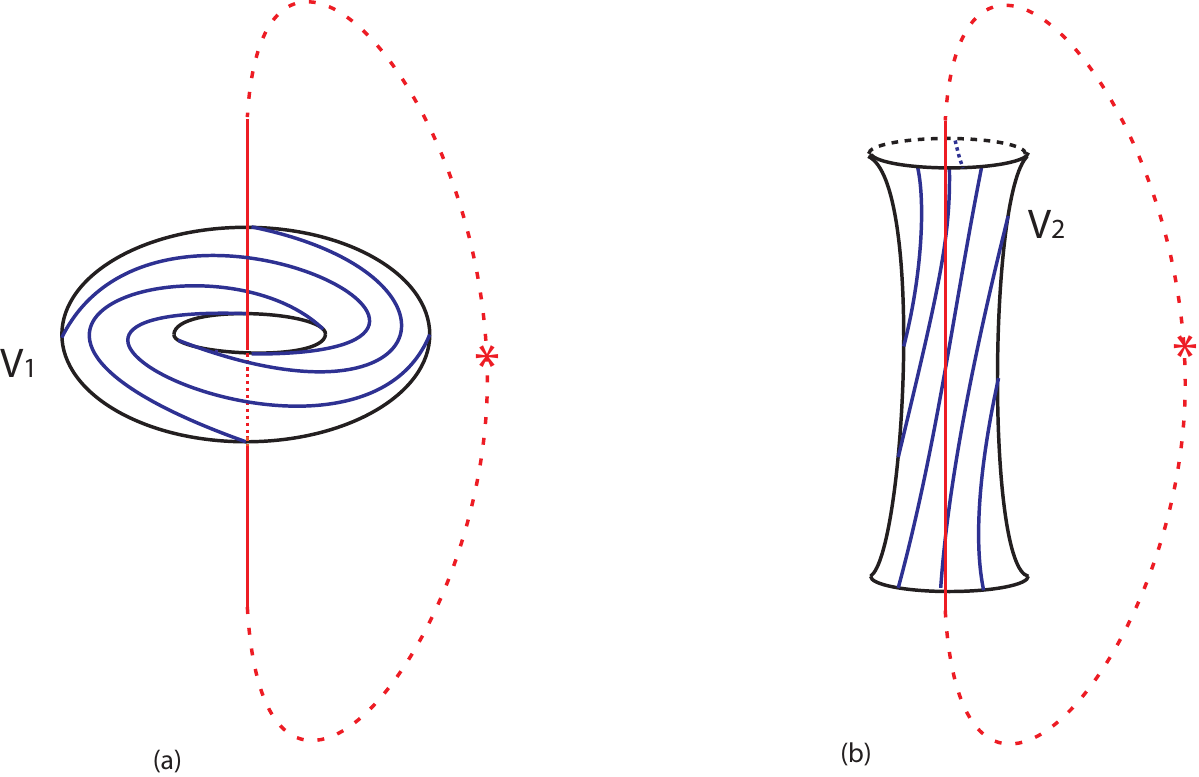}
\caption{The specified longitude becomes a meridian in the new 3-manifold.}
\label{3Dtruncatedinitial}
\end{center}
\end{figure}

Starting with $S^3$ and its description as the splitting of two solid tori, {\it 3-dimensional topological surgery} means that a solid torus $V_2 = S^1\times D^2$ is removed from $S^3$ and in the closure of the remaining manifold is replaced by another solid torus $D^2\times S^1$ (with the factors reversed), which gets attached via a homeomorphism $\phi$ along the boundary $S^1\times S^1$  of $V_2$. This boundary  (which is a torus) is the common boundary of $V_2$ with the complement solid torus $V_1$. Surgery starts and ends with two 3-manifolds and it may change the homeomorphism type of the initial 3-manifold. From the description above we obtain:
$$
M= \overline{S^3 \setminus (S^1\times D^2)} \cup_\phi (D^2\times S^1)
$$
The core of $V_2$ is called the {\it surgery curve}. Before surgery the meridians of $V_2$ bound discs, so they cut through the surgery curve (red line $l$ in Figure~\ref{3DinR3}). So, before surgery $V_2$ is layered via the indicated meridional discs.  The action of the gluing homeomorphism $\phi$ is determined by specifying  a $(p,q)$-torus knot on the boundary of $V_2$, which is  a {\it parallel curve} to the surgery curve in $V_2$. Figure~\ref{3Dtruncatedinitial}(a) illustrates a $(4,3)$-torus knot on the boundary of $V_1$. The solid torus  $V_2$ is represented by the red surgery curve, which is assumed to pass by the point at infinity. Note that, from the point of view of $V_2$ the above curve is a $(3,4)$-torus knot on the boundary of $V_2$ and it is illustrated in Figure~\ref{3Dtruncatedinitial}(b).  This $(p,q)$-torus knot is the image of the meridian via $\phi$, so it becomes a meridian in the new 3-manifold and therefore it now bounds a disc; while the meridians of $V_2$ that were bounding discs before they do not any more. See Figure~\ref{3DinR3}.
 This exchange of roles can be clearly seen in the blue parallel curve (left hand illustration) turning into a blue meridional disc (right hand illustration). So, after surgery, the layering of $V_2$ is via the discs bounded by the $(p,q)$-torus knots.  This is where we need the fourth dimension to visualize 3-dimensional surgery.  

Practically, before one could slide through a meridional disc in $V_2$ (and could also cross the surgery curve), while after surgery, the only way to come closer to the surgery curve is by following the parallel $(p,q)$-torus knot. Note that the new meridians can never reach the surgery curve.

\smallbreak

\begin{rem} \rm
Note  that the appearance  of the surgery line changes instantly the layering of the space from spheres to tori and initiates the instant creation of a simple closed curve $c$, which is the core of the solid torus $V_1$.   
\end{rem}

\begin{rem}\label{duality} \rm  
There is an apparent  duality and a natural exchange of roles of the two solid tori. Therefore, the core curve of $V_1$ could be equally considered as the surgery curve.
\end{rem}

The above topological process is called {\it $p/q$-rational surgery along the unknot} and starting from $S^3$ it results in the lens space $L(p,q)$. In fact, by a fundamental theorem of topology, every c.c.o. 3-manifold can be created from  $S^3$ by performing surgery along a knot or link (see \cite{PS, Ro}).

\subsection{}

\begin{figure}[ht!]
\begin{center}
\includegraphics[width=14cm]{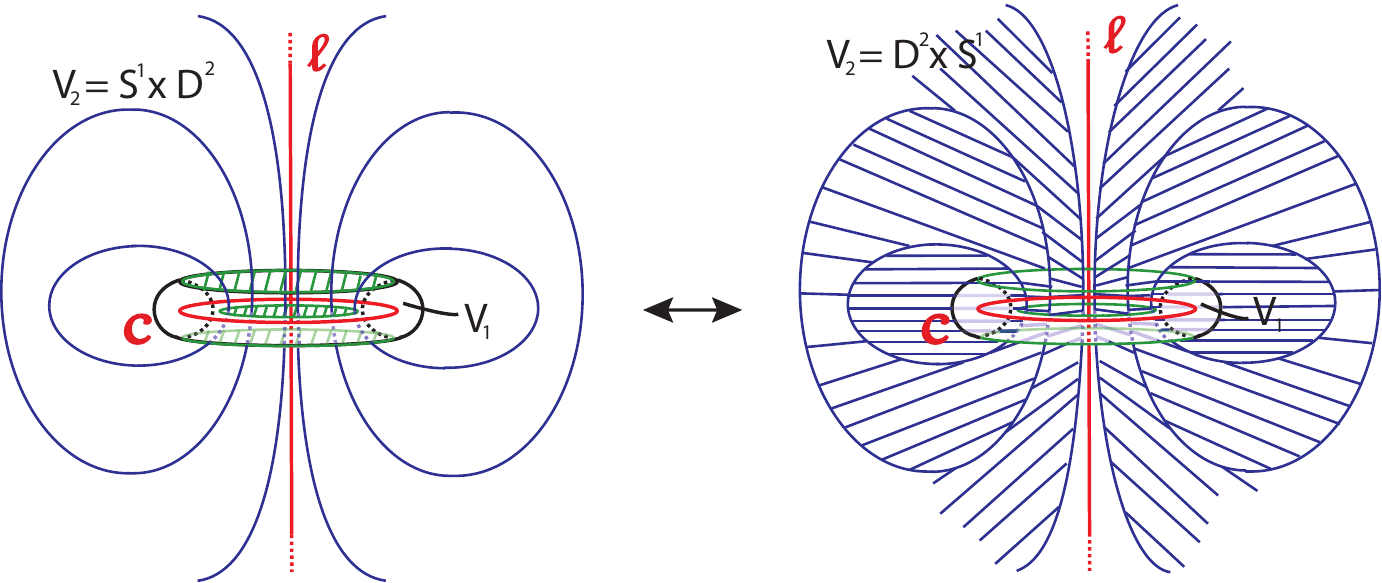}
\caption{Topological surgery along $l$.}
\label{3DinR3}
\end{center}
\end{figure}

3-dimensional surgery is much harder to visualize than lower-dimensional surgeries. A first step in this direction is to use the de-compactification of $S^3$. So, we define {\it topological surgery in} ${\mathbb R}^3$.  The only difference from the definition of surgery in  $S^3$ is that the surgery curve  is now an infinite line $l$. Figure~\ref{3DinR3} illustrates surgery in ${\mathbb R}^3$. Note that this figure resembles very much the electromagnetic field excited by a current loop which is located in the innermost torus in the drawing. Here there is no apparent drilling, that is, no obvious specified longitude, but by Remark~\ref{duality} the surgery curve is the core of the solid torus $V_1$. 

\smallbreak

\subsection{}

A second step toward visualizing 3-dimensional surgery is achieved by removing a whole  neighbourhood $B^3$ of the point at infinity. By Remark~\ref{truncate} we are then left with a 3-ball $D^3$, which is viewed as the solid torus $V_1$ corked by a (bounded) solid cylinder whose core is the arc $L$, which is part of the surgery curve. A surgery in $S^3$ along an unknotted curve passing by the point at infinity would correspond to surgery in $D^3$ along the arc $L$, the core of a solid cylinder. This solid cylinder is complemented by the outer ball $B^3$, which is homeomorphic to another solid cylinder, to form together the solid torus $V_2$. 
 The above lead to the following `localized' definition of 3-dimensional surgery. 

\begin{defn} \rm  A {\it truncated 3-dimensional surgery} in a 3-manifold $M$ is a 3-dimensional surgery, such that the surgery curve passes through the point at infinity, and such that  a neighbourhood of the point at infinity is removed. 
\end{defn}

\begin{figure}[ht!]
\begin{center}
\includegraphics[width=14cm]{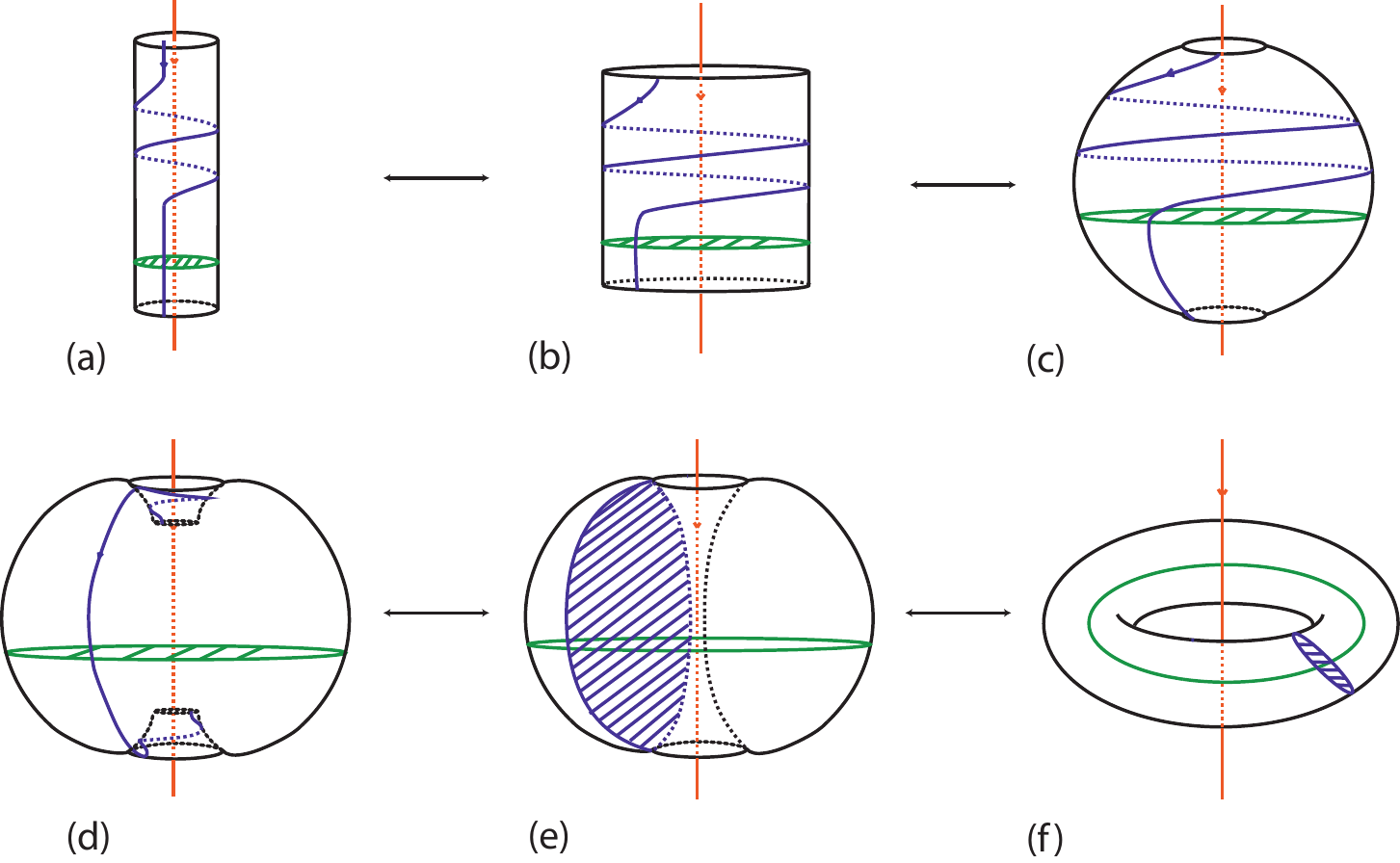}
\caption{Truncated 3-dimensional surgery helps visualize a longitude before bounding a disc afterwards.}
\label{3Dtruncated}
\end{center}
\end{figure}

This definition can help us visualize step-by-step 3-dimensional surgery along the unknot in $S^3$, especially the formation of the new meridian in the resulting lens space. For this we shall consider for simplicity a $(2,1)$-torus knot as the specified parallel curve. View Figure~\ref{3Dtruncated}. We start with a solid cylinder, which is a part  of the  solid torus $V_2$. On its boundary a $(2,1)$-curve (blue) is specified which is parallel to the core curve (red). Then the solid cylinder gets thicker and it is transformed into a 3-ball.  Then opposite parts of the cylinder move inwardly and at the same time a twisting takes place that results in `straightening' of the parallel curve. Then merging and recupling takes place resulting in a hole; thus the solid cylinder is turned into a solid torus on which the blue curve bounds now a disc. Note that the solid torus $V_1$ surrounding the initial solid cylinder is omitted in the figure.

\smallbreak

\subsection{}

Considering now life-like situations, we will proceed with inserting dynamics in truncated 3-dimensional surgery. This can help understand the topological mechanism behind some natural or physical phenomena. We start with the attracting type. 

\begin{figure}[ht!]
\begin{center}
\includegraphics[width=15cm]{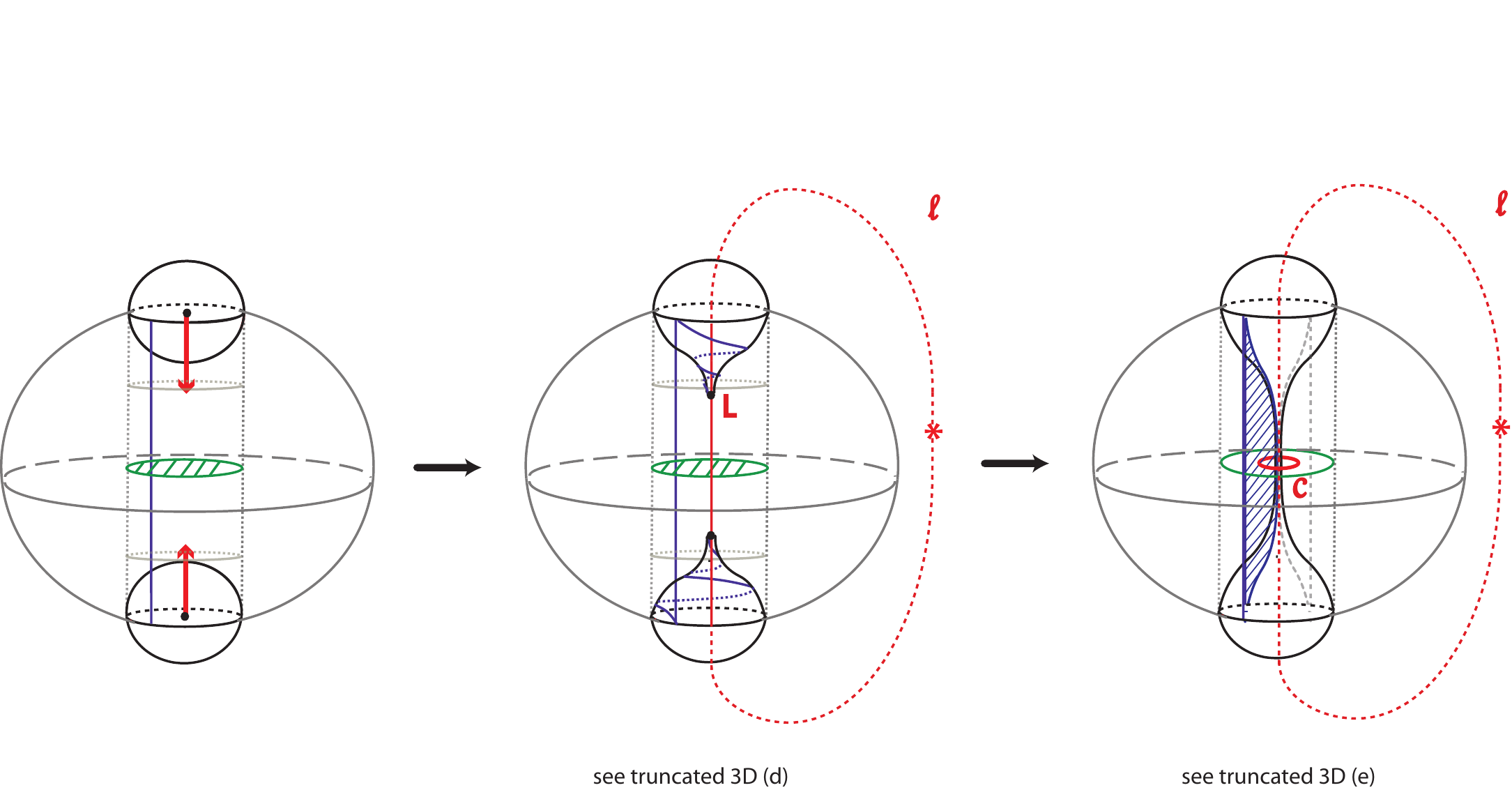}
\caption{Attracting 3-dimensional surgery.}
\label{3Dattract}
\end{center}
\end{figure}

\begin{defn} \rm 
 Consider two points in 3-space, surrounded by spherical neighbourhoods, say $B_1$ and $B_2$and assume that on these points  strong  attracting forces act. View Figure~\ref{3Dattract}.  As a result, a `joining thread', say $L$, is created between the two points and  `drilling' along $L$ is initiated. The joining arc $L$ is seen as part of a simple closed curve $l$ passing by the point at infinity. This is the surgery curve. Further,  the two 3-balls $B_1$ and $B_2$ together with the space in between make up a solid cylinder, the `cork' (cf. Figure~\ref{S3BallsToTori}). Let $V_1$ be a solid torus, which  filled by the cork gives rise to a 3-ball $D^3$, such that the centers of the two balls $B_1$ and $B_2$ lie on its boundary (compare with Figure~\ref{S3BallsToTori}). The process of attracting  3-dimensional surgery restricted in $D^3$ shall be called {\it attracting truncated 3-dimensional surgery}.

Note that the cork in the above definition is complemented by a solid cylinder, a tubular neighbourhood $B^3$ of the arc $l-L$, to the solid torus $V_2$, the complement of $V_1$ in $S^3$. This completes our familiar picture. We shall then define {\it repelling truncated 3-dimensional surgery} to be the dual phenomenon to the above, whereby strong repelling forces are applied on the two points, so strong as to initiate attracting surgery in the complementary 3-ball $B^3$, along the segment  $l-L$ with central point the point at infinity.
\end{defn}

\smallbreak

\subsection{}

\begin{figure}[ht!]
\begin{center}
\includegraphics[width=8cm]{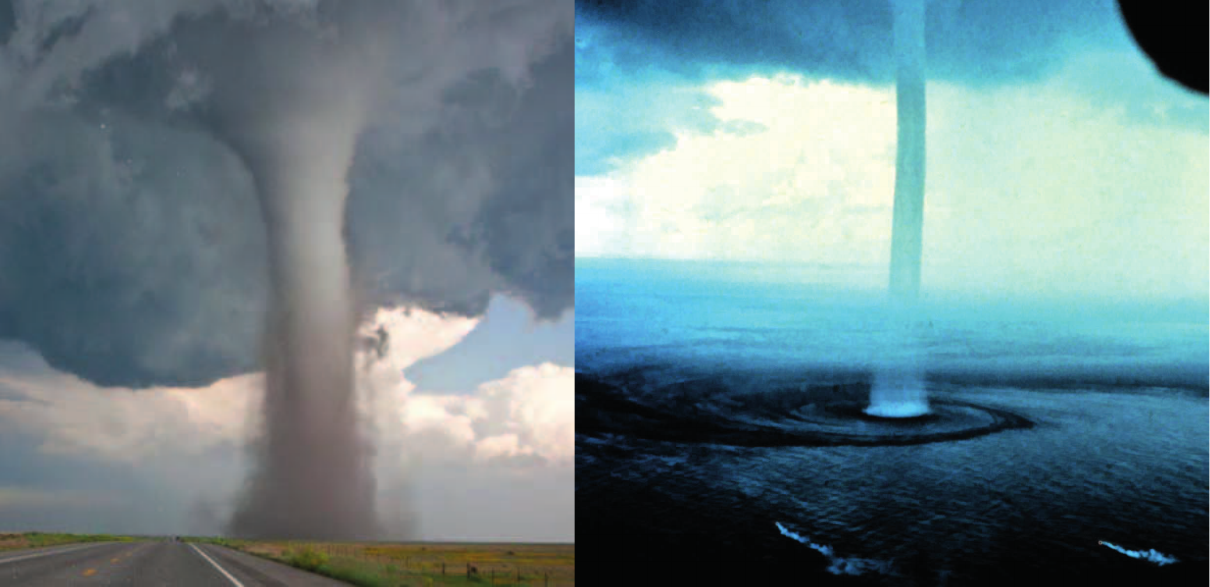}
\caption{Funnel clouds drilling and tornado formation.}
{\footnotesize \underline{{\it Sources}}: http://www.smartsuburbansurvival.com/category/natural-disasters/tornadoes.dpbs and NOAA (http://www.photolib.noaa.gov/htmls/wea00308.htm)}
\label{tornado}
\end{center}
\end{figure}

Structural similarities exhibited on vastly different scales of the universe can be visualized and explained with  3-dimensional surgery. A natural phenomenon resembling strongly the  process of truncated 3-dimensional surgery is the formation of tornadoes, see Figure~\ref{tornado}. Namely, if certain meteorological conditions are met, funnel-shaped clouds start descending toward the ground. Once they reach it, they become tornadoes. Drawing the analogy to 3-dimensional surgery, first the poles are chosen, one on the tip of the cloud and the other on the ground, and they seem to be joined through an invisible line. Then, starting from the first point, the wind revolves in a  helicoidal motion toward the second point, resembling `hole drilling' along the line until the hole is drilled. Topologically speaking, in this case seems to be undergoing rational surgery along the unknot. 

\smallbreak

\begin{figure}[ht!]
\begin{center}
\includegraphics[width=6.5cm]{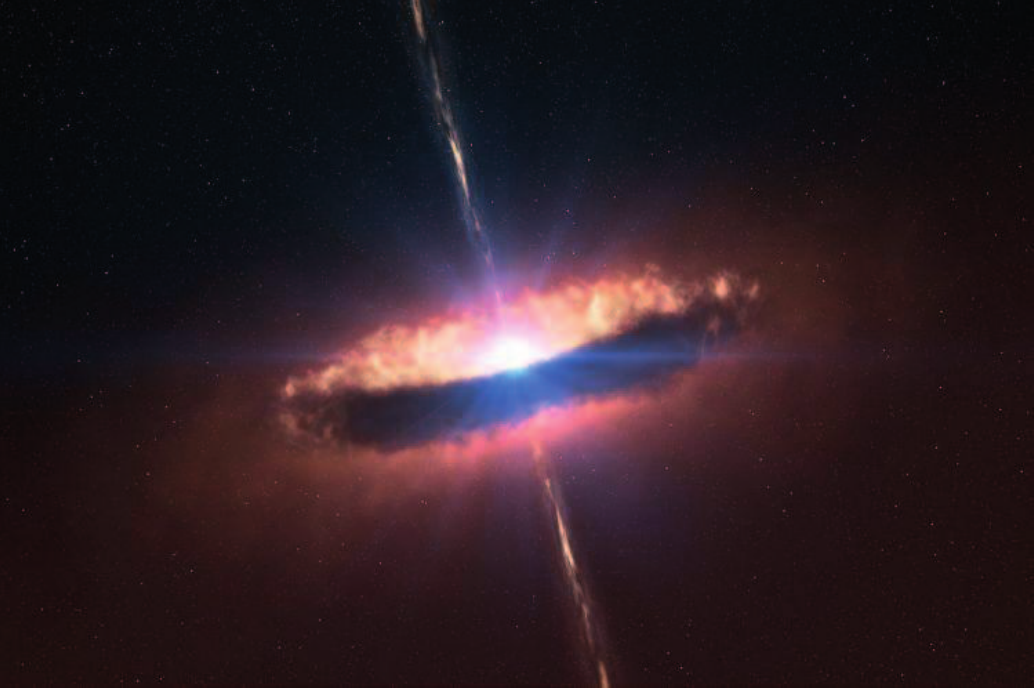}
\caption{Birth of baby star.}
{\footnotesize \underline{{\it Source}}: http://www.spitzer.caltech.edu/news/1153-feature10-11-Unravelling-the-Mystery-of-Star-Birth-Dust-Disk-Discovered-Around-Massive-Star}
\label{BabyStar}
\end{center}
\end{figure}

 There are  other examples exhibiting topological behaviour of 3-dimensional surgery. Figure \ref{BabyStar}, for example, illustrates ``a dusty disc closely encircling a massive baby star".

\smallbreak

\section{Connecting surgery in different dimensions}

Note that solid 2-dimensional surgery can be almost viewed as the intermediate stage of 3-dimensional  surgery. Indeed, there is a great resemblance between solid 2-dimensional surgery and truncated 3-dimensional  surgery. They both begin with a solid ball and there is `drilling' occuring along a `cork', a solid cylinder passing through the center. In fact, by Definition~\ref{continuum2D}, the solid 2-dimensional surgery is responsible for the creation of the curve $c$ in  truncated 3-dimensional  surgery. Yet,  there is a crucial difference:  in solid 2-dimensional surgery the cylindrical cork is removed afterwards and we left with just the solid torus $V_1$ with its core curve $c$. While, in truncated 3-dimensional  surgery, matter is still there (surrounding the arc $L$)  but it is completely altered.  The above descriptions explain the connection between 2-dimensional and 3-dimensional topological surgery, up to now not so explicitly presented. The meeting ground is the three-space with solid 2-dimensional surgery on the one end and truncated 3-dimensional surgery on the other end. 

\begin{figure}[ht!]
\begin{center}
\includegraphics[width=13cm]{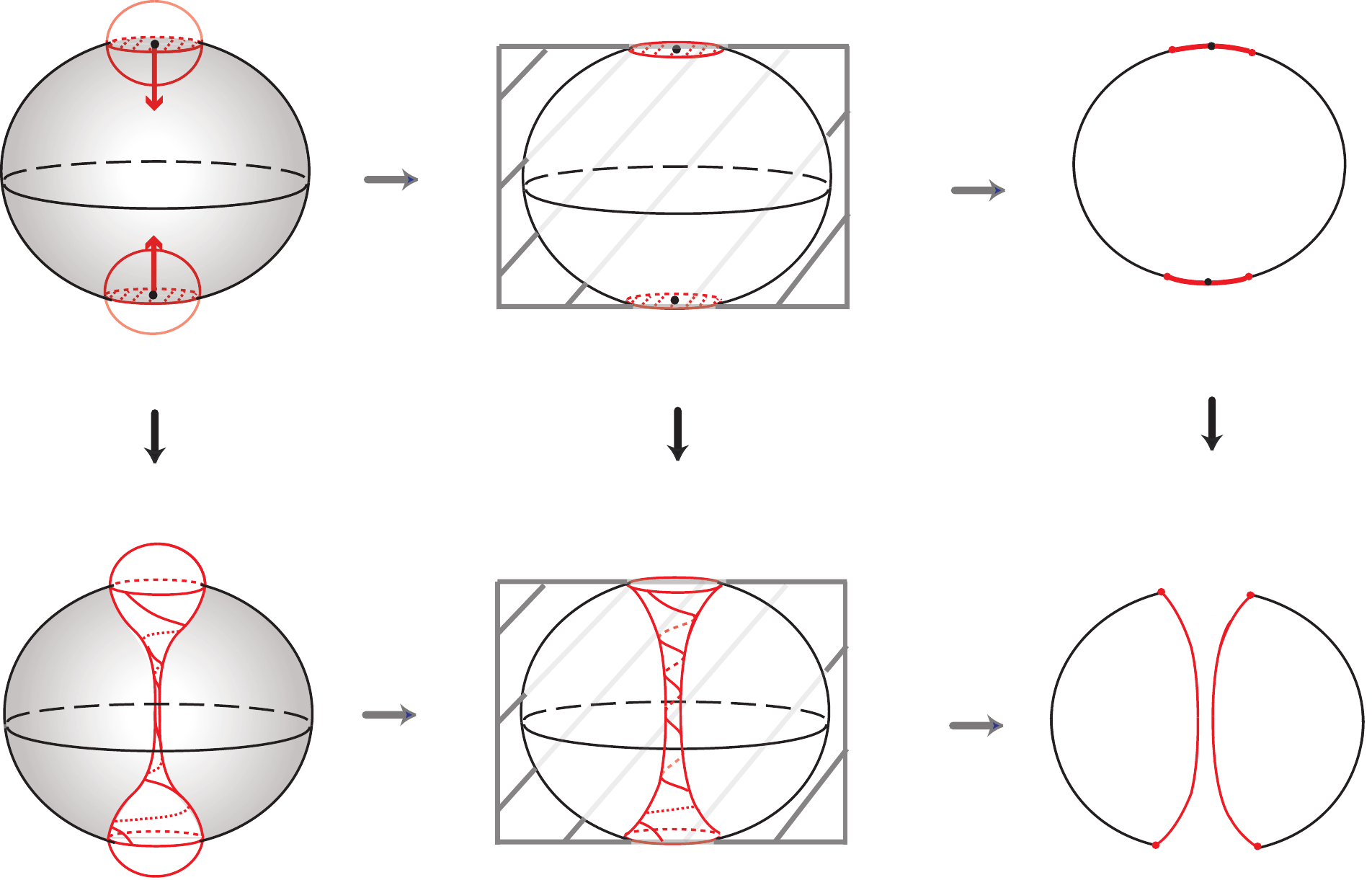}
\caption{Connecting low-dimensional surgeries.}
\label{crossections}
\end{center}
\end{figure}

We shall now go a bit further than that and explain the connection of attracting surgeries in all three dimensions. View Figure~\ref{crossections}. On the left-hand top and bottom pictures we see truncated 3-dimensional surgery. Taking on the top picture the intersection with the boundary of the 3-ball $D^3$ we pass to the initial picture of attracting 2-dimensional surgery, where two points with surrounding discs are specified. Restricting truncated 3-dimensional surgery only to this sphere results in the final stage of  attracting 2-dimensional surgery (middle bottom illustration). Taking finally the intersection with a meridional plane gives rise to 1-dimensional surgery (rightmost illustrations).

\smallbreak

\section{Conclusions}

Topological surgery occurs in numerous natural phenomena of varying scales where two points (poles) are selected and attracting or reppeling forces are applied. Examples of such phenomena comprise: DNA recombination, magnetic reconnection, mitosis, gene transfer, the creation of Falaco solitons, the formation of whirls and tornadoes and magnetic fields. 

In this paper we tried to pinn down the connection of such phenomena with topological surgery. In order to do this we first enhanced the static description of topological sugery  of dimensions~1, 2 and 3 by introducing dynamics by means of attracting or repelling forces between two `poles'. We then filled in the interior space in 1- and 2-dimensional  surgery, introducing the notion of solid 1- and 2-dimensional surgery. This way more natural phenomena can be accommodated in the connection. Finally we fitted many more natural phenomena in the connection by  intoducing the notion of truncated 1-, 2-, and 3-dimensional topological surgery, whereby surgery is more localized. Thus,  instead of considering surgery as an abstract topological process, it can now be viewed as a property of many natural phenomena. 

On the other hand, all these new notions enabled us understand and visualize 3-dimensional surgey and reveal the relation between topological surgeries in all three lower  dimensions.  In \cite{SSN} these notions are used for connecting 3-dimensional topological surgery with a dynamical system. Then, phenomena related to 3-dimensional surgery could be modelled by this dynamical system.

We hope that through this study, topology and dynamics of natural  phenomena as well as  topological surgery may now be better understood and that our connections will serve as ground for many more insightful observations.

\end{document}